\documentclass[10pt, a4paper]{article}

\usepackage{fullpage}
\usepackage{amssymb}
\usepackage{amsmath}
\usepackage{amsfonts}
\usepackage{subcaption}
\usepackage{multirow}
\usepackage{tikz}
\usepackage{float}

\newtheorem{remark}{Remark}
\newtheorem{proposition}{Proposition}

\usepackage{xcolor}
\usepackage[colorlinks]{hyperref}
\definecolor{darkgreen}{rgb}{0.0,0.4,0.0}
\definecolor{darkred}{rgb}{0.6,0.0,0.0}
\definecolor{darkblue}{rgb}{0.0,0.0,0.5}
\definecolor{gray}{rgb}{0.5,0.5,0.5}
\definecolor{cyan}{rgb}{0.0,1.0,1.0}
\definecolor{darkcyan}{rgb}{0.0,0.5,0.5}
\definecolor{darkorange}{rgb}{0.8,0.4,0.0}
\definecolor{darkmargenta}{rgb}{0.5,0.0,0.5}
\definecolor{black}{rgb}{0.0,0.0,0.0}

\hypersetup{citecolor=blue, linkcolor=red, anchorcolor=darkblue,
filecolor=darkblue, urlcolor=darkblue}
\hypersetup{pdfauthor={Nicolae Suciu},pdftitle={Suciu}}

\begin{document}

\title{Space-time upscaling of reactive transport in porous media (extended version)}

\author{Nicolae Suciu$^{1}\footnote{Corresponding author. \textit{Email adress: suciu@math.fau.de}}$ \hspace{0.1cm}, Florin A. Radu$^2$, Iuliu S. Pop$^3$\\ \\
\smallskip 
\textit{$^{1}$Tiberiu Popoviciu Institute of Numerical Analysis, Romanian Academy,}\\
\textit{Fantanele 57, 400320 Cluj-Napoca, Romania}\\
\textit{$^{2}$ Center for Modeling of Coupled Subsurface Dynamics, University of Bergen,}\\
\textit{ All\'{e}gaten 41, 5007 Bergen, Norway}\\
\textit{$^{3}$ Data Science Institute, Hasselt University,}\\
\textit{ Campus Diepenbeek, 3590 Diepenbeek, Belgium}}

\date{}
\maketitle

\begin{abstract}
Reactive transport in saturated/unsaturated porous media is numerically upscaled to the space-time scale of a hypothetical measurement through coarse-grained space-time (CGST) averages. The reactive transport is modeled at the fine-grained Darcy scale by the actual number of molecules involved in reactions which undergo advective and diffusive movements described by global random walk (GRW) simulations. The CGST averages verify identities similar to a local balance equation which allow us to derive expressions for the flow velocity and the intrinsic diffusion coefficient in terms of averaged microscopic quantities. The latter are further used to verify the CGST-GRW numerical approach. The upscaling approach is applied to biodegradation processes in saturated aquifers and variably saturated soils and the CGST averages are compared to classical volume averages. One finds that if the process is characterized by slow variations in time, as in homogeneous reaction systems, the differences between the two averages are negligible. Instead, the differences are significant and can be extremely large in simulations of time-dependent biodegradation processes in both soils and saturated aquifers. 
\end{abstract}

\begin{flushleft}

\end{flushleft}

\section{Introduction}
\label{sec:intro}

Concentrations measured in soil experiments as well as in samples pumped from monitoring wells or collected by multilevel samplers are necessarily space averages over the support scale of the sampling procedure and measuring device \cite{Maillouxetal2003,Teutschetal2000,Bayer-Raichetal2004}.
The sampling volume enters as an important parameter in analytical and numerical models because of its impact on the  moments and the
probability distribution of the uncertain concentration transported in heterogeneous natural porous media \cite{Andricevic1998,Schwedeetal2008,ToninaandBelin2008,Srzicetal2013,BarrosandFiori2014}.

The spatial average over the sampling volume is accounted for in different ways depending on modeling approaches. In analytical and
semianalytical approaches, the average concentration is an \textit{integral of a continuous field} of concentrations
\cite{Andricevic1998,Schwedeetal2008} or related probability densities \cite{BarrosandFiori2014}. The same procedure is applied
to both resident and flux-weighted concentrations \cite{Schwedeetal2008}. When the purpose is to estimate probability
distributions, the resident concentration is preferred \cite{Schwedeetal2008,BarrosandFiori2014}. The resident concentrations solve local balance equations for reactive transport which are the starting point in deriving evolution equations for the concentration probability density function \cite{Suciuetal2015} and the corresponding concentration moments equations \cite{Schueleretal2016}. These equations describe the statistics of the concentration without considering a sampling volume, which corresponds to a point sampling in volume averaging approaches \cite{ToninaandBelin2008,Srzicetal2013}. However, by averaging the solution of the local balance equation with a spatial filter, a filtered density function approach is obtained, which explicitly takes into account the sampling volume, determined by the width of the spatial filter \cite{Suciuetal2016}.

In the case of discrete descriptions of the transport process, such as sequential particle tracking \cite{CaroniandFiorotto2005,Srzicetal2013,ToninaandBelin2008,Wrightetal2021} and global random walk (GRW) simulations \cite{Suciuetal2016,Suciu2019,Suciuetal2021,SuciuandRadu2021}, concentrations are measured by \textit{counting particles} in the sampling volume. When using moderate numbers of particles, as in particle tracking approaches, the traditional counting particles produces non-smooth fields. This issue is overcome in kernel density estimations \cite{Fernandez-GarciaandSanchez-Vila2011,Sole-MariandFernandez-Garcia2018,Sole-Marietal2019} by weighting the contributions of the individual particles with smooth kernel functions. The resulting kernel estimates are smooth concentration fields and concentration derivatives/gradients. The parameter of the kernel function is a measure of the support volume which defines the region over which a particle influences the estimate. The optimal volume support can be evaluated as function of the number of particles and it has been shown that it goes to zero for infinite number of particles \cite{Fernandez-GarciaandSanchez-Vila2011}. The need of an optimal support volume is thus related to the relatively small numbers of particles used in particle tracking approaches, which is limited by the computational resources. In GRW methods which, instead of tracking individual particles, move groups of particles over the sites of a regular lattice, particles can be associated to molecules and a mole of solute is represented by Avogrado's number \cite{SuciuandRadu2021}.

The relation between the two spatial averaging procedures, by integrating continuous fields and by counting particles, is not trivial. In fact, this is equivalent to passing from a microscopic level of description to a macroscopic one. Such connections can be rigorously established if the spatial sampling is replaced by a space-time average \cite{Suciu2001,Vamos2007}.

Considering the time scale along with the support volume is imposed by the obvious fact that any measurement is also a time average and hydrological measurements are no exception \cite{DestouniandGraham1997,Andricevic1998}. Even if some flow and transport conditions (e.g. small spatial scale of
concentration fluctuations \cite{Andricevic1998} or steady state transport \cite{Schwedeetal2008}) justify a sampling volume approach which disregards the measurement time scale, in most cases both the spatial and temporal scales specific to the observation method \cite{Bayer-Raichetal2004,Berkowitz2021,DestouniandGraham1997} have to be considered.

The influence of the observation method on the expected value and the variance of the locally measured concentrations has been modeled in \cite{DestouniandGraham1997} by time averages of Lagrangian concentration representations. The proposed methodology can be used for either unsaturated or saturated flow conditions and for different types of observation methods. The temporal scale may be related to a spatial scale, as for instance for constant water volume sampling or for constant soil sample length. But the temporal scale may also be independent of the spatial scale, as in case of monitoring wells, where it is related to the pumping time \cite{Bayer-Raichetal2004}, or in lysimeter experiments \cite{DabrowskaRykala2021}, where the relevant scale is the effluent collection time. The analysis based on temporal averages presented in \cite{DestouniandGraham1997} does not consider the effect of the spatial averaging and the results are consistent with those reported in previously published volume sampling studies.

A theoretical approach based on spatial and temporal averages was developed to transform microscale balance equations into macroscopic equations \cite{HeandSykes1996}. The basic assumption in this approach is the existence of differential equations which describe the transport of the volumetric density of an extensive physical quantity at the microscale. The averaging domain consists of the set product between a representative elementary volume which may vary in space and time and a variable timescale. The space-time averaging domain can thus be thought of as a generalization of the
classical ``fluid molecule''. Assuming the existence of the solution of the microscale equation, averaging theorems which relate derivatives of averages to averages of derivatives are proved and further used to derive the macroscopic equation for multiphase transport in porous media \cite{Cushman1983,HeandSykes1996}. This equation is not closed and constitutive relations are needed to formulate well-posed problems.

Most recent papers on averaging procedures for porous media do not consider the temporal averaging, focusing mainly on volume averaging theorems for multiphase systems (see \cite{GrayandMiller2013} and references therein). Spatial averaging procedures are used, for instance, to formulate macroscopic equations for multiphase transport \cite{Grayetal2015}, to understand how the internal structure of a two-fluid porous medium system influences the
macroscale state \cite{McClure2017}, to account for dependencies between the size of the averaging volume and macroscale in derivations of Darcy's law from Navier-Stokes equations \cite{Nordbottenetal2007}, as well as in upscaling transport with Monod reactions to obtain effective equations for the macroscale \cite{Hesseetal2009}.

Important achievements in deriving upscaled models for porous media were obtained by using homogenization techniques \cite{Hornung1996} which,
unlike direct volume or spatio-temporal averaging, lead to closed upscaled equations with macroscopic coefficients completely determined by solving unit-cell problems \cite{AuriaultandAdler1995}. Examples encompass but are not limited to, derivation of the Darcy scale model for precipitation and dissolution in porous media \cite{Kumaretal2016}, modeling injection of cold water in geothermal reservoirs \cite{Bringedaletal2016}, upscaling flow by homogenization methods for non-periodic and random parameters of the subsurface system \cite{Nolenetal2008}. The concept of spatio-temporal upscaling in the context of homogenization by multiple-scale expansions, initially developed to study linear equations with coefficients oscillating in both space and time \cite{BensoussanLionsPapanicolaou2011}, has been extended to nonlinear diffusion equations \cite{AkagiandOka2020} and to reactive transport in porous media driven by oscillating pore water velocity \cite{VanDuijnandvanderZee2018}. In this approach, spatial and temporal scales are necessarily related, which could be a limitation of the spatio-temporal homogenization from the point of view of engineering applications. Another inconvenience is that the spatial and temporal scales are related to the inverse of the small parameter of the homogenization problem and the macroscopic description obtained when this parameter goes to zero cannot be directly related to the scales of measurement instruments and methods \cite{DestouniandGraham1997,Bayer-Raichetal2004} or to the spatial and temporal scales of the hydrological processes \cite{Skoienetal2003}.

In both direct spatio-temporal averaging and homogenization approaches the ``microscopic'' state is given by continuous fields. The alternative ``counting particles'' averaging approach requires a microscopic description given by trajectories of physical or computational particles. Assuming that the heterogeneity of the subsurface system can be modeled by random functions, a ``microscopic'' description corresponding to the local, Darcy scale of the transport process is given by trajectories of diffusion in random velocity fields \cite{Suciu2019} modeled numerically by particle tracking \cite{Wrightetal2021} or GRW algorithms \cite{Suciuetal2021}. A general method to obtain a macroscopic description at spatio-temporal scales that can be directly related to those of the measurements and hydrological observations is the coarse-grained (CG) description of corpuscular systems through space-time (ST) averages \cite{Vamosetal1996}. The CGST approach also yields unclosed macroscopic equations but, since a complete microscopic description is available, it provides a frame to investigate the structure of the constitutive relations through numerical simulations \cite{Vamosetal1997}.

Upscaling reactive transport in subsurface hydrological systems aims at obtaining models that can account for small scale processes in predicting the fate of contaminants at larger scales and in designing remediation strategies for contaminated sites. Solutions proposed so far, based on either Eulerian approaches \cite{Portaetal2012} or on particle tracking transport models combined with Eulerian descriptions of the reactions \cite{Wrightetal2021} are essentially volume averaging methods. However, it has been argued that small-scale temporal fluctuations of reacting species concentrations, in conditions of locally varying effects of advection and diffusion, propagate from local to larger spatial scales and at least an implicit temporal averaging must also be included in upscaling reactive transport \cite{Berkowitz2021}. In the cited paper, it is further argued that the temporal information can be included through a continuous-time random walk approach, with a power-law distribution of transition times calibrated from measurements, which accounts for temporal effects over a range of time scales.

The CGST method offers the perspective of explicitly relating the upscaled concentrations to the spatio-temporal scales of interest. The fine-grained concentrations are represented microscopically by the actual number of molecules involved in reactions, their advective-diffusive movement being described by simple discrete-time random walk processes modeled with GRW algorithms \cite{SuciuandRadu2021}. This molecular description implicitly accounts for spatio-temporal fluctuations produced by inhomogeneities of reaction and transport parameters. Moreover, by using integrated flow and transport GRW solvers for unsaturated/saturated porous media \cite{Suciuetal2021} the CGST approach can be applied to upscaling contaminant transport in both soil and aquifer systems.

In this article we propose a CGST upscaling approach using microscopic descriptions consisting of one- and two-dimensional GRW simulations of the reactive transport at Darcy scale. General statistical mechanics results on CGST averaging are presented in Section~\ref{sec:cgst_av}. In Section~\ref{sec:verification}, theoretical formulas for the flow velocity and the intrinsic diffusion coefficient are derived as corollaries of the general result and are used to verify the Matlab codes designed for the numerical space-time upscaling. Further, in Section~\ref{sec:STupscaling} the CGST approach is applied to various one-dimensional problems of reactive transport solved by GRW algorithms, with emphasis on biodegradation processes in saturated and unsaturated natural porous media. Two-dimensional codes for CGST upscaling are verified and applied to biodegradation in variably saturated soils in Section~\ref{sec:cgst_2dim}. The main results and summarized in Section~\ref{sec:concl}. The proof of Proposition 1 presented in Section~\ref{sec:cgst_av} is given in Appendix~A and Appendix~B illustrates the influence of the space-time scale on the discrepancy between volume and CGST averages. The codes implemented in Matlab used in this study are stored in the Git repository \href{https://github.com/PMFlow/SpaceTimeUpscaling}{SpaceTimeUpscaling} \cite{Suciuetal2023}.

\section{Coarse-grained space-time averages}
\label{sec:cgst_av}

The CGST method allows continuous modeling of corpuscular systems, without the need to use the vague notion of ``fluid particle'' which is small enough to be infinitesimal with respect to the size of the macroscopic system, but large enough to contain many molecules in local thermodynamic equilibrium. Continuous macroscopic fields and balance equations can be derived under the only assumption that a system of particles (e.g. molecules and other physical particles, or computational particles used in GRW and particle tracking simulations) possesses a kinematic microscopic description through piecewise analytic time functions. A general proof of this result has been obtained in \cite{Vamosetal1996} by considering CGST averages over spheres and symmetric time intervals. Proofs for one-dimensional systems were presented in \cite{Vamosetal1997,Vamosetal2000} and a general result for averaging on hypercubes in the phase space of the statistical mechanics, for conserved number of phase points, was achieved in \cite{Vamos2007}. Considering a kinematic description for systems with a variable number of particles, we prove in Appendix~A the smoothness properties of the CGST averages over $d$-dimensional cubes and symmetric time intervals stated in the following proposition.

\begin{proposition}
\label{proposition:cg_average}
If a discrete system of $\mathcal{N}$ particles is described by piecewise analytic time functions $\varphi_{i}(t):[0,T]\longmapsto\mathbb{R}$,
$i=1,\dots,\mathcal{N}$, representing trajectory and velocity components, $x_{\alpha i},\; \xi_{\alpha i}=dx_{\alpha i}/dt,\; \alpha=1,\ldots,d$, or other physical properties associated to the $i$-th particle, then

1) there exists a macroscopic description
of the same system given by almost everywhere (a.e.) continuous
fields through \textit{coarse-grained space-time averages}
\begin{equation}\label{eq:average}
\langle\varphi\rangle(\mathbf{x},t;a,\tau)=\frac{1}{2\tau(2a)^d}\sum_{i=1}^{\mathcal{N}}\int_{t-\tau}^{t+\tau}\varphi_i(t')\prod_{\alpha=1}^{d}
\chi_{\alpha}(x_{\alpha i}(t'))dt',
\end{equation}
where $\tau<T/2$ is the half length of the temporal averaging interval, $(2a)^d$ is the volume of the open $d$-dimensional cube $C(\mathbf{x},a)=\{\mathbf{x},\mathbf{y}\in\mathbb{R}^d\; | \; |y_{\alpha}-x_{\alpha}|<a, \alpha=1,2,\ldots,d\}$, and $\chi_{\alpha}$ is the characteristic function of the interval $(x_{\alpha}-a,x_{\alpha}+a)$,

2) $\langle\varphi\rangle(\mathbf{x},t;a,\tau)$ has continuous partial derivatives
a.e. in $\mathbb{R}^{d}\times(\tau,T-\tau)$ and satisfies the
identity
\begin{equation}\label{eq:cg_eq}
\partial_{t}\langle\varphi\rangle+\nabla\cdot\langle\varphi\boldsymbol{\xi}\rangle=\textstyle{\left\langle\frac{d\varphi}{dt}\right\rangle} +\delta\varphi,
\end{equation}
where $\delta\varphi$ accounts for discontinuous variations produced when a particle is created or consumed in chemical reactions or by finite jumps of the function $\varphi_{i}(t)$.
\end{proposition}

For $\varphi _{i}\equiv1$, $i=1\ldots\mathcal{N}$, the CGST average (\ref{eq:average}) counts the number of particles inside the cube $C(\mathbf{x},a)$ during a time interval of length $2\tau$, that is, $\langle 1\rangle$ is an a.e. continuous approximation of the macroscopic concentration field. The average over the number of particles  defined by $\overline{\boldsymbol{\xi}}=\langle\boldsymbol{\xi}\rangle/\langle 1\rangle$, if $\langle 1\rangle>0$, and by $\overline{\boldsymbol{\xi}}=0$, if $\langle 1\rangle=0$ provides an a.e. continuous approximation of the macroscopic velocity field \cite{Vamosetal1996,Vamosetal2000}. In particular, if the particles move with the constant velocity $\mathbf{u}$, according to (\ref{eq:average}), $\overline{\boldsymbol{\xi}}=\mathbf{u}$. In the absence of chemical reactions and jump discontinuities the term $\delta\varphi$ in (\ref{eq:cg_eq}) vanishes and one obtains a
continuity equation,
\begin{equation}\label{eq:conti}
\partial_{t}\langle 1\rangle+\nabla\cdot(\langle 1\rangle\overline{\boldsymbol{\xi}})=0.
\end{equation}

For $\tau\longrightarrow 0$ the CGST average (\ref{eq:average}) with $\varphi _{i}\equiv1$ gives the estimation of the concentration by
the usual volume averaging formula. The essential advantage of (\ref{eq:average}) is that the additional time average yields a.e. continuous fields
obeying the identity  (\ref{eq:cg_eq}), which provides a frame to close the balance equations \cite{Vamosetal1997,Vamosetal2000}.

In the framework of the statistical mechanical theory of transport processes \cite{IrvingandKirkwood1950}, the continuous concentration field $c(\mathbf{x},t)$ at the scale of the mathematical continuum, hereafter called ``fine-grained concentration'', is defined as
an average $M$ over the ensemble of realizations of the transport process,
\[
c(\mathbf{x},t):=M\left[\sum\limits_{i=1}^{\mathcal{N}}\delta(\mathbf{X}_{i}(t,\omega)-\mathbf{x})\right],
\]
where $\delta$ is the Dirac function and $\mathbf{X}_{i}(t,\omega)$ denotes the realization $\omega$ of the random trajectory. Applying the averaging operator $M$ to (\ref{eq:average}), for $\varphi_{i}\equiv1$, with $x_{\alpha i}(t)=X_{\alpha i}(t,\omega)$, and expressing the characteristic function of the cube in (\ref{eq:average}) with the aid of the Dirac functional
\[
\prod_{\alpha=1}^{d}\chi_{\alpha}(X_{\alpha i}(t))=\chi_{C(\mathbf{x},a)}(\mathbf{X}_i,t)=
\int_{\mathbb{R}^d}\chi_{C(\mathbf{x},a)}(x')\delta(\mathbf{X}_i(t,\omega)-\mathbf{x}')d\mathbf{x}',
\]
one obtains the following result.

\begin{proposition}
\label{proposition:stoch_average}
The ensemble average of the coarse-grained average $\langle 1\rangle(\mathbf{x},t;a,\tau)$ is an
average  of the fine-grained concentration $c(\mathbf{x},t)$ over the cube
$C(\mathbf{x},a)$ and the time interval $(t-\tau,t+\tau)$ given by
\begin{equation}\label{eq:stoch_average}
M[\langle 1\rangle](\mathbf{x},t;a,\tau)=\frac{1}{2\tau(2a)^d}\int\limits_{t-\tau}^{t+\tau}dt^{\prime}\int\limits_{C(\mathbf{x},a)}c(\mathbf{x}',t')d\mathbf{x}'.
\end{equation}
\end{proposition}
The ensemble average $M[\langle 1\rangle](\mathbf{x},t;a,\tau)$ defines the continuous concentration field $c(\mathbf{x},t;a,\tau)$ as observed at the spatio-temporal scale $(a,\tau)$. Equation (\ref{eq:stoch_average}) is a particular case of the general relation
between the smooth continuous field and the CGST average of its microscopic description \cite[Eq.~(4.6)]{Suciu2001}. For $\tau\longrightarrow 0$,
Eq. (\ref{eq:stoch_average}) gives the relation between the ``integrating continuous fields'' and ``counting particles'' averaging procedures in sampling volume approaches.

The use of CGST averages (\ref{eq:average}) in numerical simulations is conditioned by the piecewise analyticity of $\varphi (t)$. Random walks on lattices are eligible because the corresponding trajectories $\mathbf{x}_{i}(t)$ and velocities $\boldsymbol{\xi}_{i}$ are piecewise analytic functions with numerable discontinuity points at the jump instants \cite{Vamosetal2000}. In GRW simulations, the sum over the number of particles in (\ref{eq:average}) is equivalent to summing up the contributions of the groups of particles jumping between the same lattice sites in the interior of the cube $C(\mathbf{x},a)$, in the time interval $(t-\tau,t+\tau)$.

\section{Verification of the CGST averaging procedure}
\label{sec:verification}

The goal of the space-time upscaling will be achieved by performing CGST averages on microscopic descriptions consisting of simulations of transport in porous media performed with GRW algorithms. In this section, we consider the one-dimensional case and we verify the reliability of the averaging procedure applied to simulations of diffusion with both the unbiased GRW and its biased version, BGRW.

\subsection{Implications of the CGST upscaling}
\label{sec:upscaling}
In the one-dimensional case, the CGST average (\ref{eq:average}) becomes
\begin{equation}\label{eq:average1}
\langle\varphi\rangle(x,t;a,\tau)\;=\;\frac{1}{4\tau a}
\;\sum\limits_{i=1}^{\mathcal{N}}\;\int\limits_{t-\tau}^{t+\tau}\,\varphi_{i}(t')\,\chi(x_{i}(t'))\,dt'
\end{equation}
and the continuity equation (\ref{eq:conti}) takes the form
\begin{equation}\label{eq:conti1}
\partial_{t}\langle 1\rangle+\partial_x(\langle 1\rangle\overline{\xi})=0.
\end{equation}
Considering $\varphi_i=x_i$ and using (\ref{eq:cg_eq}), the advective flux in (\ref{eq:conti1}) can be expressed as
\[
\langle 1\rangle\overline{\xi}=\langle\xi\rangle=\textstyle{\left\langle\frac{dx}{dt}\right\rangle}=\partial_{t}\left\langle x\right\rangle+\partial_{x}\left\langle x\xi\right\rangle.
\]
Further, by assuming $\langle 1\rangle>0$ and by considering the average $\overline{x}=\langle x\rangle/\langle 1\rangle$ of particle positions in the averaging domain $(x-a, x+a)\times(t-\tau, t+\tau)$ and the continuity equation (\ref{eq:conti1}), we get
\[
\partial_{t}\left\langle x\right\rangle=\partial_{t}(\langle 1\rangle\overline{x})=\langle 1\rangle\partial_t\overline{x}+\langle 1\rangle\overline{\xi}\partial_x\overline{x}-\partial_x(\langle 1\rangle\overline{x}\overline{\xi})
\]
and the continuity equation takes the form of Fokker-Planck equation
\begin{equation}\label{eq:adv_diff}
\partial_{t}\langle 1\rangle+\partial_{x}[\langle 1\rangle(\partial_{t}\overline{x}+\overline{\xi}\partial_x\overline{x})]=\partial^{2}_{x}[(\overline{x}\overline{\xi}-\overline{x\xi})\langle 1\rangle].
\end{equation}
The positivity of the coefficient $(\overline{x}\overline{\xi}-\overline{x\xi})$ in (\ref{eq:adv_diff}) is a ``microscopic criterion of irreversibility'' for the thermodynamic system of the $\mathcal{N}$ particles \cite{Suciu2001}. However, the drift coefficient $(\partial_{t}\overline{x}+\overline{x}\partial_x\overline{x})$ differs from the Eulerian drift $u$ and Eq.~(\ref{eq:adv_diff}) cannot be generally identified with the advection-diffusion equation valid at the macroscopic scale. Instead, assuming $\partial_{t}\langle x\rangle=0$, and $u=0$, which implies $\partial_{t}c=0$ according to (\ref{eq:conti1}), Eq.~(\ref{eq:adv_diff}) reduces to
\begin{equation}\label{eq:diff}
\partial^{2}_{x}(D\langle 1\rangle)=0,
\end{equation}
which is the stationary diffusion equation with constant diffusion coefficient $D=-\overline{x\xi}=-\left\langle x\xi\right\rangle/\langle 1\rangle$.

\begin{remark}
\label{rem:instrinsicD}
If the random walkers are relatively homogeneously distributed inside the spatial averaging interval $(x-a,x+a)$, then $\overline{x}\approx x$, implying $\partial_t \overline{x}=0$ and $\partial_x \overline{x}=1$. With these, Eq.~(\ref{eq:adv_diff}) becomes
\[
\partial_{t}\langle 1\rangle-\partial_x(\langle 1\rangle\overline{\xi})-x\partial^{2}_{x}(\langle 1\rangle\overline{\xi})=\partial^{2}_{x}(D\langle 1\rangle),
\]
that is, an advection-diffusion equation with coefficient $D$, introduced in Eq.~(\ref{eq:diff}) above, but with an unusual drift term, which is again different from that in the advection-diffusion equation known from the axiomatic formulation of the fluid dynamics, i.e.~$\partial_{t}c+\partial_x(cu)=\partial^{2}_{x}(Dc)$. Actually, this advection-diffusion equation is just another reformulation, under the hypothesis $\overline{x}\approx x$, of the continuity equation (\ref{eq:conti1}) verified by the CGST averaged number of particles $\langle 1\rangle$. The merit of these reformulations is that they highlight the self-diffusion process with ``intrinsic diffusion coefficient'' $D$ of the system of particles of the microscopic description.
\end{remark}
The relation defining the intrinsic diffusion coefficient $D=-\left\langle x\xi\right\rangle/\langle 1\rangle$, as well as the average defining the macroscopic velocity, $u=\left\langle \xi\right\rangle/\langle 1\rangle$, will be used in the following to verify the numerical approach for computing CGST averages over microscopic descriptions consisting of systems of random walkers.

\subsection{GRW/BGRW algorithms}
\label{sec:GRWalg}

The GRW algorithms can be understood as weak schemes for It\^{o} stochastic differential equations \cite{Suciu2019}. In the above one-dimensional case, a weak scheme for solving the equation $dX_t=udt+\sqrt{2D}dW$, where $W$ is the standard Wiener process, can be obtained by using the discrete random walk process $X_{k\Delta t}$ with increments $X_{(k+1)\Delta t}-X_{k\Delta t}=(\tilde{u}+\zeta d)\Delta x$ consisting of jumps over $(\tilde{u}+\zeta d)$ steps $\Delta x$ on a regular lattice. The integer number $d$ is the amplitude of the diffusive displacements, $\tilde{u}=\lfloor u\Delta t/\Delta x\rfloor$ is the truncated advective displacement, $\zeta$ is a random variable taking three discrete values with probabilities $P(\zeta=\pm 1)=r/2$, $P(\zeta=0)=(1-r)$, and the real parameter $r=2D\Delta t/(d\Delta x)^2$ is the dimensionless diffusion coefficient. The probability density of the It\^{o} process $X_t$, i.e. the normalized concentration of random walkers $c(x,t)$ solving the Fokker-Planck equation $\partial_t c + \partial_x(uc)=D\partial^{2}_{x}c$, is approximated by counting the number of random walkers at lattice sites $l$, $c(l\Delta x,k\Delta t)\thickapprox n_{l,k}/\mathcal{N}$, where $\mathcal{N}$ is the total number of random walkers. As in Section~\ref{sec:cgst_av} above, by $c(x,t)$ we denote the fine-grained concentration.

The unbiased GRW algorithm distributes the random walkers over lattice sites according to
\begin{equation}\label{eq:grw}
n_{l,k}=\delta n_{l+\tilde{u}|l,k}+\delta n_{l+\tilde{u}-d|l,k}+\delta n_{l+\tilde{u}+d|l,k},
\end{equation}
where the numbers $\delta n$ of random walkers are binomial random variables with statistics specified by the parameter $r$. They represent the number of particles which at time $k$ undergo jumps from a lattice site $l$ to $l+\tilde{u}$, $l+\tilde{u}-d$, and $l+\tilde{u}+d$, respectively.

An alternative weak scheme for the above It\^{o} equation uses the random walk process with increments $X_{(k+1)\Delta t}-X_{k\Delta t}=\zeta\Delta x$ and biased jump probabilities $P(\zeta=\pm 1)=(r\pm u)/2$, $P(\zeta=0)=(1-r)$. This results in a biased algorithm, BGRW, which distributes the random walkers on first-neighbor lattice sites according to
\begin{equation}\label{eq:bgrw}
n_{l,k}=\delta n_{l|l,k}+\delta n_{l-1|l,k}+\delta n_{l+1|l,k}.
\end{equation}
The restriction $|u|\le r$ implies an upper bound for the local P\'{e}clet number, $\text{P\'{e}}=|u|\Delta x/D\le 2$, which imposes smaller $\Delta x$ and $\Delta t$, rendering the BGRW algorithm slower than the unbiased GRW algorithm.

For a general presentation of the GRW/BGRW algorithms we refer to \cite[Chap. 3]{Suciu2019}. Implementation details, for more complex problems consisting of one-component reactive transport fully coupled to saturated/unsaturated flows governed by the degenerate Richards equation can be found in \cite[Sect. 4.2]{Suciuetal2021}, and for transport with multi-species Monod reactions in \cite[Sect. 3]{SuciuandRadu2021}.

\subsection{Averaging procedure and verification}
\label{sec:GRWaverage}

The CGST average (\ref{eq:average1}) on the GRW/BGRW simulations is computed by summing up contributions from groups $\delta n$ of random walkers as follows. Let $N_a$ be the number of lattice sites inside the open interval $(x-a,x+a)$. If the conditions $|X_{k\Delta t}-x|\leq a$ and $|X_{(k+1)\Delta t}-x|<a$ are fulfilled at every occupied site, the CGST average of the microscopic quantity $\varphi$ is computed by
\begin{equation}\label{eq:average1_grw}
\langle\varphi\rangle(x,t;a,\tau)\;=\;\frac{\Delta t}{4\tau a}
\;\sum\limits_{k=(t-\tau)/\Delta t}^{(t+\tau)/\Delta t}\,\sum\limits_{i=1}^{N_{a}}\;\Phi_{i,k}(\varphi),
\end{equation}
where $\Phi_{i,k}(\varphi)$ are contributions of the quantity $\varphi$ to the CGST average.

As microscopic quantities we consider in the following $\varphi_k=1$, for the computation of the CGST concentration $\langle 1\rangle$, $\varphi_k=X_{(k+1)\Delta t}$, for the mean particle positions $\langle x\rangle$, and $\varphi_k=(X_{(k+1)\Delta t}-X_{k\Delta t})/\Delta t$, for the mean velocity $\langle \xi\rangle$. The computation of the intrinsic diffusion coefficient $D=-\langle x\xi\rangle/\langle 1\rangle$ defined during the derivation of the stationary diffusion equation (\ref{eq:diff}) will be also used in the following to verify the correctness of the averaging procedure. Since $D$ is a property of the thermodynamic system (Remark~\ref{rem:instrinsicD}), it does not depend on the flow regime and can be computed as well for microscopic descriptions consisting of GRW simulations of advection-diffusion processes. To do that, we have to consider the velocity $\xi$ of the random walk jumps. In the case of the microscopic description provided by the unbiased GRW algorithm, $\xi$ computed by disregarding the advective displacement $\tilde{u}$ is associated to the middle of the jump interval and $\langle x\xi\rangle$ is computed for the microscopic quantity $\varphi_k=(X_{k\Delta t}+\tilde{u}\Delta x+\frac{1}{2}\zeta d\Delta x)(\zeta d\Delta x)/\Delta t$. When using the BGRW algorithm (\ref{eq:bgrw}), where the bias of the random walk jump probabilities accounts for advection, $\langle x\xi\rangle$ is computed with binomial random variables $\delta n_{l\pm 1\mid l,k}$ representing the number of random walkers jumping on neighboring sites  with unbiased probabilities $\pm r/2$, where $r$ is the same parameter as that used to compute biased jump probabilities $(r\pm u)/2$ in case of non-vanishing advection velocity $u$, and the corresponding microscopic quantity is given by $\varphi_k=(X_{k\Delta t}+\frac{1}{2}\zeta\Delta x)(\zeta \Delta x/\Delta t)$.

Using (\ref{eq:grw}) and (\ref{eq:bgrw}), the CGST averages considered in this article are computed according to (\ref{eq:average1_grw}) with the functions $\Phi_{i,k}(\varphi)$ listed below:
\begin{itemize}
\item the concentration $\langle 1\rangle$

$\Phi_{i,k}(\varphi)=n_{i,k+1}$
\item the mean particle positions $\langle x\rangle$
\begin{itemize}
\item GRW algorithm: $\Phi_{i,k}(\varphi)=[(i+\tilde{u})\Delta x]n_{i,k}$
\item BGRW algorithm: $\Phi_{i,k}(\varphi)=(i\Delta x)\delta n_{i|i,k}+[(i-1)\Delta x]\delta n_{i-1|i,k}+[(i+1)\Delta x]\delta n_{i+1|i,k}$
\end{itemize}
\item the mean velocity $\langle \xi\rangle$
\begin{itemize}
\item GRW algorithm: $\Phi_{i,k}(\varphi)=[\tilde{u}(\Delta x/\Delta t)]n_{i,k}$
\item BGRW algorithm: $\Phi_{i,k}(\varphi)=(-\Delta x/\Delta t)\delta n_{i-1|i,k}+(\Delta x/\Delta t)\delta n_{i+1|i,k}$
\end{itemize}
\item the product position-velocity $\langle x\xi\rangle$
\begin{itemize}
\item GRW algorithm:

$\Phi_{i,k}(\varphi)=[(i+\tilde{u}-\frac{1}{2}d)\Delta x](-d\Delta x/\Delta t)\delta n_{i+\tilde{u}-d|i,k}+[(i+\tilde{u}+\frac{1}{2}d)\Delta x](d\Delta x/\Delta t)\delta n_{i+\tilde{u}+d|i,k}$
\item BGRW algorithm:

$\Phi_{i,k}(\varphi)=(i-\frac{1}{2})\Delta x(-\Delta x/\Delta t)\delta n_{i-1|i,k}^{(unbiased)}+(i+\frac{1}{2})\Delta x(\Delta x/\Delta t)\delta n_{i+1|i,k}^{(unbiased)}$

\end{itemize}
\end{itemize}

The averaging procedure is verified by computing the flow velocity and the intrinsic diffusion coefficient with the expressions derived by CGST averaging in Section~\ref{sec:upscaling}. We shall consider the fine scale transport process described by the advection-diffusion equation with constant coefficients $D$ and $u$,
\begin{equation}\label{eq:fine_grainedAD}
\partial_t c + u\partial_x(c)=D\partial^{2}_{x}c.
\end{equation}
The fine-grained concentration $c$ solving (\ref{eq:fine_grainedAD}) is approximated numerically with biased and unbiased GRW algorithms, which also provide microscopic descriptions of the physical system. As follows from Eq.~(\ref{eq:average1}), in this case the macroscopic flow velocity equals the constant drift $u$. Also, one expects that the intrinsic diffusion coefficient equals the diffusion parameter $D$ of the advection-diffusion process.

For illustration and verification purposes, we consider the one-dimensional spatial domain $\Omega=(0,1)$, a time interval of length $T=1$, and GRW/BGRW lattices of constant $\Delta x$ and $\mathcal{L}=1/\Delta x +1$ sites, denoted by $l=1,\ldots, \mathcal{L}$. We impose no-flux conditions $n_{1,k}=n_{2,k}$ and $n_{\mathcal{L},k}=n_{\mathcal{L}-1,k}$ which ensure a steady-state solution of Eq.~(\ref{eq:fine_grainedAD}). The initial condition consists of $\mathcal{N}=10^{24}$ random walkers uniformly distributed on the lattice, that is, a constant concentration $n_{l,k}/\mathcal{N}=1/\mathcal{L}$ mole per lattice site. Under these conditions, according to (\ref{eq:average1_grw}), one obtains the normalized CGST average concentration $\langle 1\rangle/\mathcal{N}=N_a/(2a\mathcal{L})$, which does not depend on $\tau$ and equals the spatial average over the interval $(x-a,x+a)$. With $a=0.05$ and $\tau=0.1$, the CGST averages are computed over 9 disjoint spatial intervals of length $2a$ and three time intervals of length $2\tau$ centered at the sampling times $t=$0.1, 0.5, and 0.9. The lattice constant is progressively decreased from $\Delta x=5\cdot 10^{-3}$ to $\Delta x=3.33\cdot 10^{-4}$ and the corresponding time steps $\Delta t$ are obtained from the condition $r=2D\Delta t/(d\Delta x)^2\leq 1$, with $d=1$ in case of BGRW algorithm, and with both $d=1$ and $\tilde{u}=1$ in case of unbiased GRW (See Section~\ref{sec:GRWalg}). The constant velocity is set to $u=1$ and the diffusion coefficient is set to $D=0.0001$.

\begin{figure}
\begin{minipage}[t]{0.45\linewidth}\centering
\includegraphics[width=\linewidth]{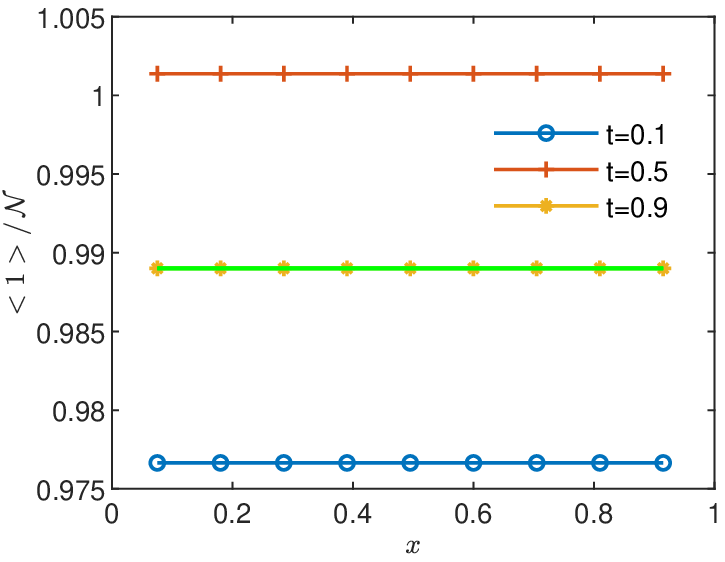}
\caption{\label{fig:ver_c1000}Normalized coarse-grained concentrations compared to the spatial average $N_a/(2a\mathcal{L})$ (full line without markers) computed with BGRW for $\Delta x=10^{-3}$.}
\end{minipage}
\hspace*{0.1in}
\begin{minipage}[t]{0.45\linewidth}\centering
\includegraphics[width=\linewidth]{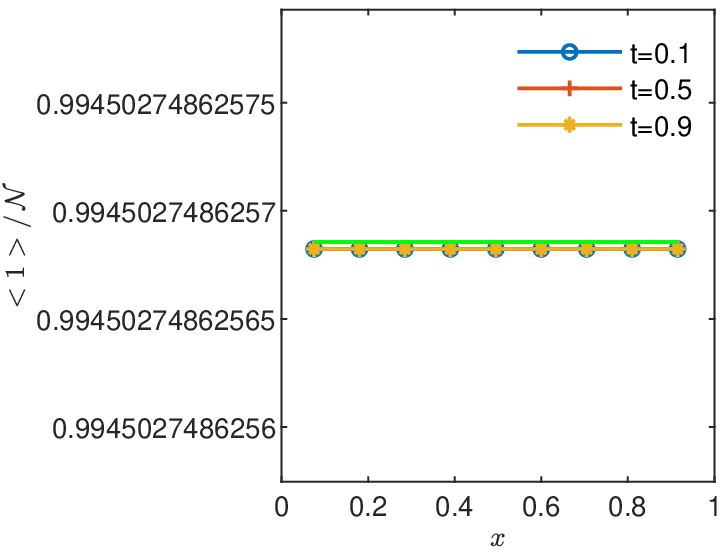}
\caption{\label{fig:ver_c2000}Normalized coarse-grained concentrations compared to the spatial average $N_a/(2a\mathcal{L})$ (full line without markers) computed with BGRW for $\Delta x=5\cdot 10^{-4}$.}
\end{minipage}
\end{figure}

\begin{figure}
\begin{minipage}[t]{0.45\linewidth}\centering
\includegraphics[width=\linewidth]{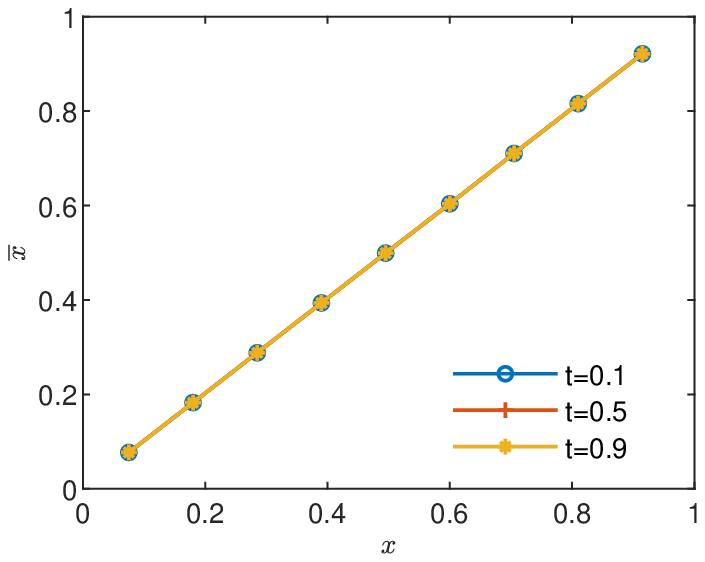}
\caption{\label{fig:ver_r1000}CGST average of particle positions computed with BGRW for $\Delta x=10^{-3}$.}
\end{minipage}
\hspace*{0.1in}
\begin{minipage}[t]{0.45\linewidth}\centering
\includegraphics[width=\linewidth]{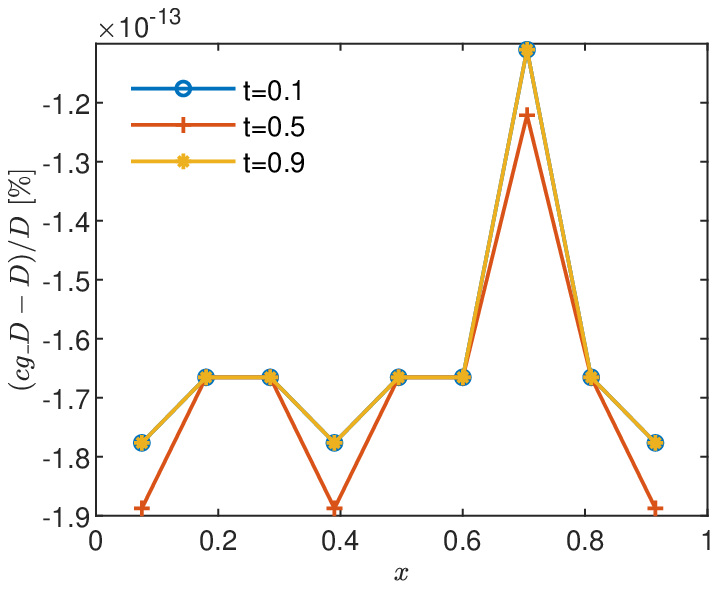}
\caption{\label{fig:ver_errD1000}Relative errors of the coarse-grained diffusion coefficient $cg\_D$ computed with BGRW for $\Delta x=10^{-3}$.}
\end{minipage}
\end{figure}

It has been found that the steady state, as indicated by a CGST concentration almost identical to that given by the spatial average, is reached for $\Delta x\leq 5\cdot 10^{-4}$ in case of the BGRW microscopic description (see Figs.~\ref{fig:ver_c1000}~and~\ref{fig:ver_c2000}) and for $\Delta x\leq 3.33\cdot 10^{-4}$ in case of the GRW description. However, as seen in Figs.~\ref{fig:ver_r1000}~and~\ref{fig:ver_errD1000}, the mean particle positions are stationary, with $\overline{x}\approx x$, and the coarse-grained computation of the intrinsic diffusion coefficient $cg\_D$ is highly accurate, even if the coarse-grained concentration is not yet stationary. The mean diffusion coefficient $\overline{D}$, the mean velocity $\overline{u}$, and the corresponding standard deviations are calculated from sequences of $cg\_D$ and $u$ values at the $N_a$ points inside the spatial averaging interval averaged over the three sampling times. The results presented in Table~\ref{table:Du} are highly accurate for both $cg\_D$ and $u$.

Since the implementation of the CGST approach with the BGRW microscopic description is simpler and free of overshooting errors in simulations of biodegradation reactions \cite[Sect. 5.2]{SuciuandRadu2021}, it will be used in the following to illustrate the space-time upscaling in reactive transport problems.

\begin{table}
\caption{\label{table:Du}CGST computation of the intrinsic diffusion coefficient and velocity.}
\begin{center}
\begin{tabular}{ c c c | c c c }
    & \hspace{2cm}$\overline{D}$  & & \hspace{2cm}$\overline{u}$ \\
\hline
    $\Delta x$ & BGRW & GRW  & BGRW & GRW \\
  \hline
     5.00e-03  & 0.0001$\pm$ 2.25e-19 & 0.0001$\pm$ 2.51e-19 & \hspace{-0.6cm}1.00$\pm$ 0.00 & 1.00$\pm$ 0.00 \\
     2.50e-03  & 0.0001$\pm$ 5.07e-19 & 0.0001$\pm$ 6.21e-19 & \hspace{-0.6cm}1.00$\pm$ 0.00 & 1.00$\pm$ 0.00 \\
     1.25e-03  & 0.0001$\pm$ 1.16e-18 & 0.0001$\pm$ 1.30e-18 & 1.00$\pm$ 2.36e-16 & 1.00$\pm$ 0.00 \\
     6.25e-04  & 0.0001$\pm$ 1.42e-18 & 0.0001$\pm$ 1.78e-18 & 1.00$\pm$ 2.36e-16 & 1.00$\pm$ 0.00 \\
  \hline
\end{tabular}
\end{center}
\end{table}

\section{CGST averaging in one-dimensional transport problems}
\label{sec:STupscaling}

We consider the one-dimensional advection-diffusion transport in saturated porous media of two reactive chemical species governed at the fine-grained scale by the coupled system of equations
\begin{equation}\label{eq:react}
\partial_t c_{\nu} + \partial_x(c_{\nu}u)-D\partial^{2}_{x}c_{\nu}=R_{\nu}(c_{1},c_{2}),\; \nu=1,2.
\end{equation}

As shown in \cite{DestouniandGraham1997}, the space-time scales of the experimental measurements are determined by the observation method. The application of the CGST approach in modeling real experiments is the subject of forthcoming studies. In the following, we illustrate the CGST space-time upscaling by hypothetical experiments with independent spatial and temporal scales which ensure a good resolution of the averages and are
bellow the space-time scales on which the reacting species concentrations undergo significant changes. Precisely, in the examples presented below, we choose the spatial scale $a$ as 1/30 of the length of the domain in the mean flow direction and the temporal scale $\tau$ as 1/11 of the total observation time.

The discrepancy between volume and CGST averages at a given time point is quantified by relative errors computed with vector $l^2$-norms as $e_{c_{\nu}}=\|\overline{c}_{\nu}-\langle 1_{\nu}\rangle\|/\|\langle 1_{\nu}\rangle\|$, as  well as by maximum relative differences
\[
\varepsilon_{c_{\nu}}=\max(|\overline{c}_{\nu}-\langle 1_{\nu}\rangle|)/\langle1_{\nu}\rangle(\arg(\max(|\langle 1_{\nu}\rangle - \overline{c}_{\nu}))), \; \nu=1,2.
\]
The influence of $a$ and $\tau$ on discrepancy is analyzed for a specific problem in Appendix B.

\subsection{Bimolecular reaction with conservation of the total mass}
\label{sec:bimolReaction}

Let us begin by considering a simple reaction where the first species is consumed and the second species is augmented at the same rate. The corresponding reaction rates are given by $R_1=-K_rc_{1}c^{*^{2}}_{2}$ and $R_2=K_rc_{1}c^{*^{2}}_{2}$, with constant $K_r$. Summing up the two equations (\ref{eq:react}) one obtains an advection-diffusion equation for the passive transport of the sum $(c_{1}+c_{2})$, that is, the total mass is conserved.

\begin{figure}
\begin{minipage}[t]{0.45\linewidth}\centering
\includegraphics[width=\linewidth]{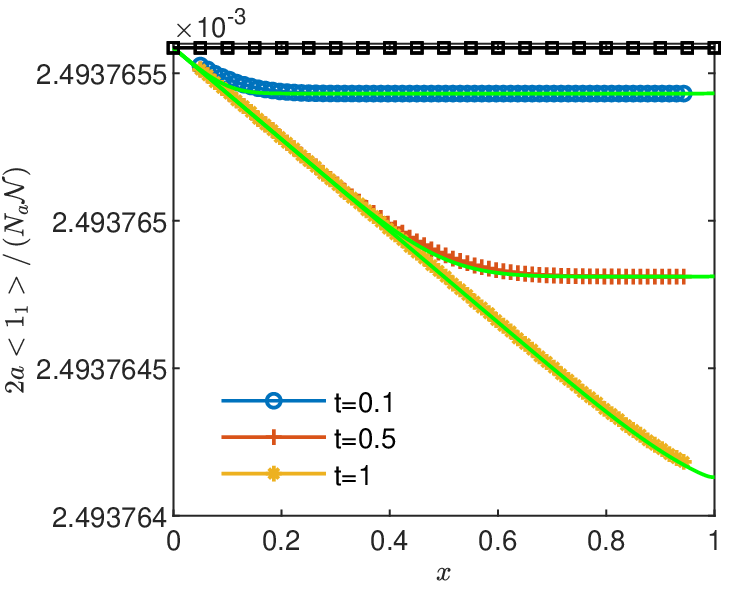}
\caption{\label{fig:bimol_c1}Moving averages, fine-grained concentration $c_{1}$ (full lines) and half total concentration $(c_{1}+c_{2})/2$ (squares).}
\end{minipage}
\hspace*{0.1in}
\begin{minipage}[t]{0.45\linewidth}\centering
\includegraphics[width=\linewidth]{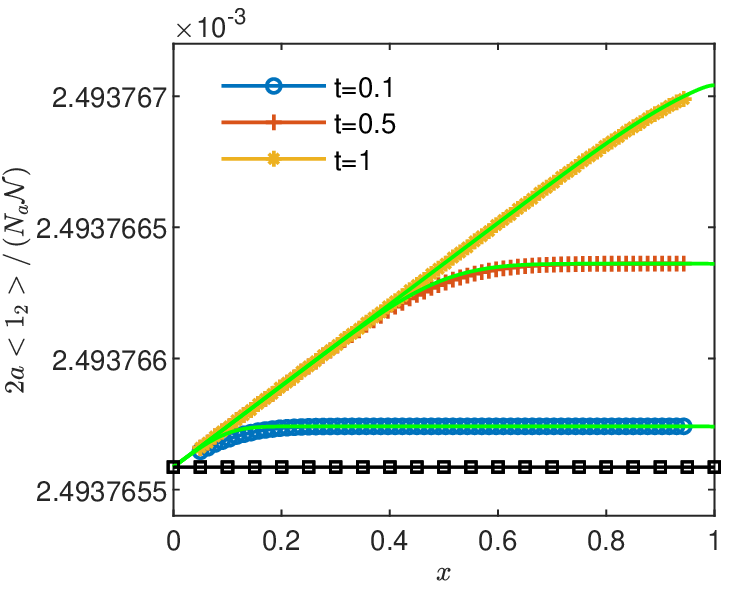}
\caption{\label{fig:bimol_c2}Moving averages, fine-grained concentration $c_{2}$ (full lines) and half total concentration $(c_{1}+c_{2})/2$ (squares).}
\end{minipage}
\end{figure}

\begin{figure}
\begin{minipage}[t]{0.45\linewidth}\centering
\includegraphics[width=\linewidth]{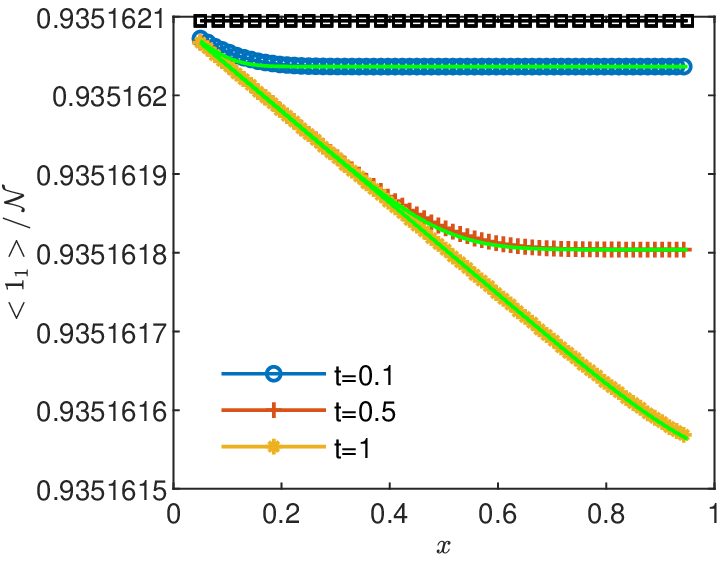}
\caption{\label{fig:bimol_CGc1}Coarse-grained concentrations of the first molecular species compared to volume averages (full lines) and half total concentration $(\langle 1_{1}\rangle+\langle 1_{2}\rangle)/2$ (squares).}
\end{minipage}
\hspace*{0.1in}
\begin{minipage}[t]{0.45\linewidth}\centering
\includegraphics[width=\linewidth]{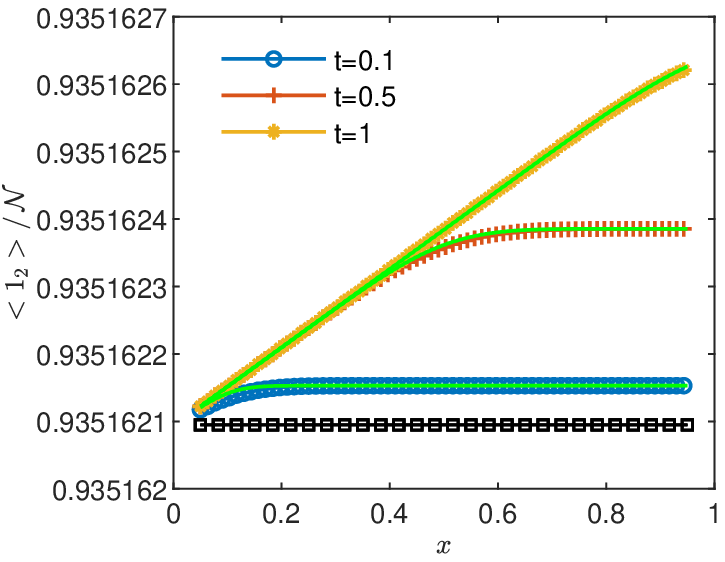}
\caption{\label{fig:bimol_CGc2}Coarse-grained concentrations of the second molecular species compared to volume averages (full lines), and half total concentration $(\langle 1_{1}\rangle+\langle 1_{2}\rangle)/2$ (squares).}
\end{minipage}
\end{figure}

Similarly to the derivation of Eq.~(\ref{eq:conti}), we find that the CGST averages of the two molecular species satisfy the particular form of Eq.~(\ref{eq:cg_eq}),
\begin{equation}\label{eq:react_cg_eq}
\partial_{t}\langle 1_{\nu}\rangle+\partial_x(\langle 1_{\nu}\rangle\overline{\xi})=\delta 1_{\nu},\; \nu=1,2,
\end{equation}
where we use the notation $1_{\nu}$ to specify the molecular species when integrating particle trajectories in Eq.~(\ref{eq:average1}) to compute the CGST average concentrations. According to (\ref{eq:average}), the terms $\delta 1_{\nu}$ account for creation and consumption of particles, hence, they can be interpreted as space-time upscaled reaction terms.

\begin{remark}\label{rem:eqCGST}
While the fine-grained concentrations $c_{\nu}$ are governed by the advection-diffusion-reaction equations~(\ref{eq:react}), the balance equations~(\ref{eq:react_cg_eq}) verified by the CGST concentrations $\langle 1_{\nu}\rangle$ have the general form of an advection-reaction equation.
\end{remark}

For consistency reasons, the conservation of the sum $(c_{1}+c_{2})$ of the fine-grained concentrations should imply the conservation of the
upscaled concentrations. The latter implies $\delta 1_{1}+\delta 1_{2}=0$ and the passive transport of the sum $(\langle 1_{1}\rangle+\langle 1_{2}\rangle)$. To verify this statement, we compute CGST averages on the BGRW description associated to Eqs.~(\ref{eq:react}), with constant parameters $u=1$, $D=0.01$, and $K_r=0.1$. The two molecular species are represented, as in Section~\ref{sec:GRWaverage} above, by one mole of particles uniformly distributed on the one-dimensional lattice at the initial time. As boundary conditions we fix the concentration on $l=1$ to its initial value and set no flux conditions on $l=\mathcal{L}$. We use the domain $\Omega=(0,1)$, total time $T=1.1$, spatial scale $a=0.03$, and temporal scale $\tau=0.1$. The reactive transport is solved by a splitting scheme which alternates BGRW transport steps with reaction steps computed deterministically (a one-dimensional version of the two-dimensional splitting approach proposed in \cite[Sect. 3.1]{SuciuandRadu2021}). The computations are carried out with $\Delta x=2.50\cdot 10^{-3}$ and $\Delta t=1.56\cdot 10^{-4}$.

\begin{table}[h]
\caption{\label{table:RelDiff_1Dbimol}Relative differences between volume and CGST averages: one-dimensional transport with bimolecular reactions.}
\begin{center}
\begin{tabular}{ c | c c | c c }
   $ t $ & $e_{c_1}$ & $\varepsilon_{c_1}$ &  $e_{c_2}$ & $\varepsilon_{c_2}$ \\
\hline
     0.1 & 2.44e-09 & 7.54e-09 & 2.44e-09 & 7.54e-09 \\
     0.5 & 1.81e-09 & 3.99e-09 & 1.81e-09 & 3.99e-09 \\
     1 & 7.54e-10 & 2.53e-09 & 7.54e-10 & 2.53e-09 \\
  \hline
\end{tabular}
\end{center}
\end{table}

The fine-grained concentrations are compared in Figs.~\ref{fig:bimol_c1}-\ref{fig:bimol_c2} with the particular centered moving average consisting of the arithmetic mean of the numbers of particles at the $N_a$ sites inside the spatial interval $(x-a,x+a)$, averaged over the time interval $(t-\tau,t+\tau)$ and normalized by the total number $\mathcal{N}$, computed with the formula $2a\langle 1_{\nu}\rangle/\,(N_a\mathcal{N})$. One finds a very good agreement between the two sets of curves, as expected for slowly varying concentrations $c_{\nu}$. 

Since the concentration is not vanishing in the domain during the computation, the intrinsic diffusion coefficient and the velocity can be calculated as in Section~\ref{sec:verification} above. One obtains accurate estimates $D=9.96\cdot 10^{-3}\pm 1.42\cdot 10^{-3}$ and $u=9.98\cdot 10^{-1}\pm 1.70\cdot 10^{-2}$, which certifies the correctness of the CGST averaging procedure. In Figs.~\ref{fig:bimol_CGc1}-\ref{fig:bimol_CGc2}, the normalized CGST averages $\langle 1_{\nu}\rangle/\,\mathcal{N}$ are compared to the spatial averages computed, similarly to (\ref{eq:average1_grw}), by $\overline{c}_{\nu}=(1/2a)\sum_{i=1}^{N_{a}}n_{i,}$ at the center $t$ of the temporal interval. In this case, as well, the agreement is very good. The conservation of the sum of the two species concentrations is verified at both the fine- and the coarse-grained scales. As shown in Table~\ref{table:RelDiff_1Dbimol}, the two averages are practically indistinguishable, with relative differences of the order of $10^{-9}$. Hence, for small temporal variations of the concentration, as in the example presented here, the volume average alone is an appropriate model of the experimental measurements.

\subsection{Reactive transport with Monod reactions}
\label{sec:MonodReaction}

In the following, we apply the CGST upscaling procedure to aerobic biodegradation processes in the subsurface \cite{Wiedemeieretal1999}. An electron donor contaminant (e.g., benzene) of concentration $c_{1}$ is consumed by the biomass contained in the aqueous system, hereafter assumed constant in time and uniformly distributed in $\Omega$ with $c_{bio}=1$, in the presence of an electron acceptor (oxygen) of concentration $c_{2}$. The electron donor and acceptor are consumed at rates $R_1=-\theta\alpha_1\mu$ and $R_2=-\theta\alpha_2\mu$, where $\theta$ is the volumetric water content, $\alpha_1$ and $\alpha_2$ are stoichiometric constants, and
\begin{equation}\label{eq:monod}
\mu=\frac{c_{1}}{M_1+c_{1}}\frac{c_{2}}{M_2+c_{2}}
\end{equation}
is a Monod term with parameters $M_1$ and $M_2$. For the purpose of the present illustration of the CGST procedure, we consider the example used in \cite{BauseandKnabner2004,Brunneretal2012} and we set the parameters of the reaction system to $\alpha_1=5$, $\alpha_2=0.5$, $M_1=M_2=0.1$, constant $\theta=1$ for aquifers, and variable $\theta$ for partially saturated soils.

\subsubsection{Monod reactions in heterogeneous aquifers}
\label{sec:MonodReaction_randV}

We consider a one-dimensional biodegradation problem in the domain $\Omega=(0,1)$ similar to the two-dimensional one proposed in \cite{Cirpkaetal1999} to investigate the role of transverse dispersion and used in \cite{BauseandKnabner2004,Klofkornetal2002,SuciuandRadu2021} to verify new numerical approaches. The reactive transport is governed by the coupled system of equations (\ref{eq:react}) with constant local dispersion coefficient $D=0.01$ and stationary velocity given by realizations of a random space function with ensemble mean $M[u]=1$. The flow velocity is approximated by the Kraichnan procedure described in \cite[Appendix C.3.2.2]{Suciu2019} as sum of 100 cosine random modes, with amplitudes,  wavenumbers and phases determined by the parameters of the random log-hydraulic conductivity field $\ln K$. The latter is a statistically homogeneous random function with isotropic Gaussian correlation, finite correlation length $\lambda=0.01$, and variance $\sigma^2=0.1$. The Kraichnan procedure is equivalent to a linearization of the flow equation for saturated aquifers and provides accurate approximations of the solution for small variances $\sigma^2$ of the statistically homogeneous $\ln K$ field \cite{Schwarzeetal2001}. The chosen number of cosine random modes ensures accurate simulations of passive transport with random velocity over traveled distances, in units of $\lambda$, of the same order of magnitude, i.e. $100\lambda$ in the present computations \cite{Eberhardetal2007}. Since the Gaussian correlation depends on the spatial variable as $\sim \exp(x/\lambda)^2$ and the saturated flow problem is linear, the transport simulation in the domain $\Omega$ of unit length also represents the solution in a domain of dimension $1/\lambda=100$ \cite{Alecsaetal2020,Suciuetal2021}. With meters and days as space and time units, this corresponds to a typical setup for contaminant transport at local scales in saturated aquifers \cite{Suciu2019}.

\begin{figure}
\begin{minipage}[t]{0.45\linewidth}\centering
\includegraphics[width=\linewidth]{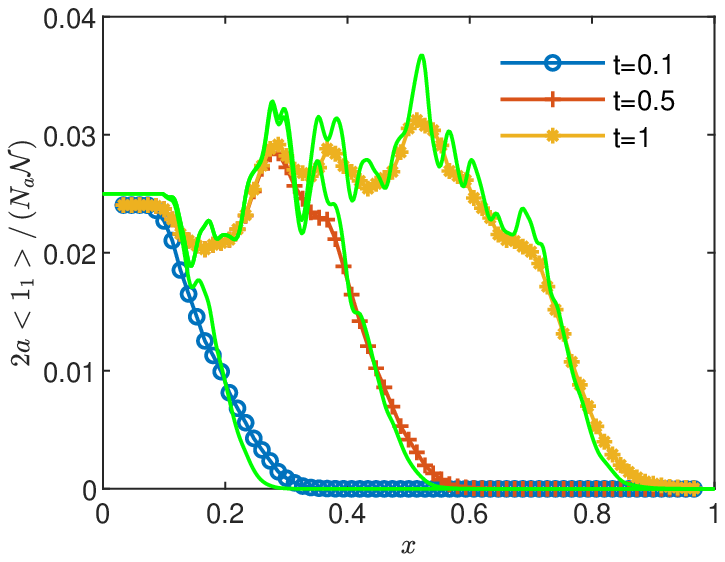}
\caption{\label{fig:rand_c1}Moving averages and fine-grained contaminant concentrations $c_{1}$ (full lines).}
\end{minipage}
\hspace*{0.1in}
\begin{minipage}[t]{0.45\linewidth}\centering
\includegraphics[width=\linewidth]{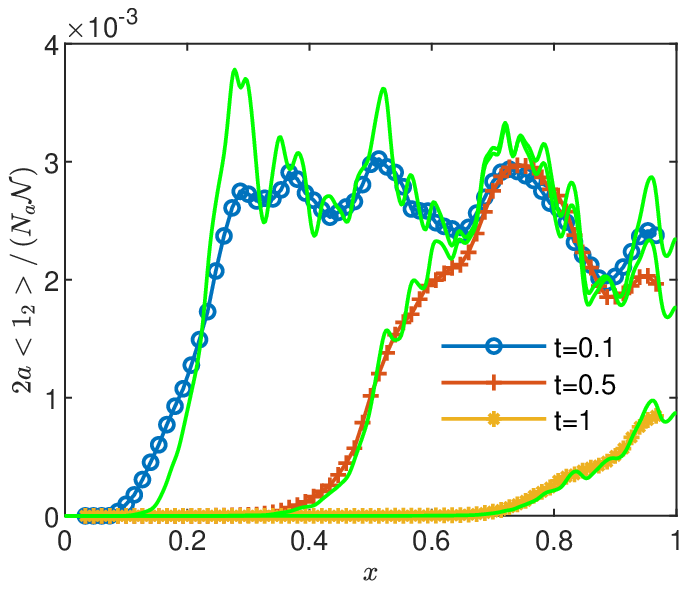}
\caption{\label{fig:rand_c2}Moving averages and fine-grained oxygen concentrations $c_{2}$ (full lines).}
\end{minipage}
\end{figure}

\begin{figure}
\begin{minipage}[t]{0.45\linewidth}\centering
\includegraphics[width=\linewidth]{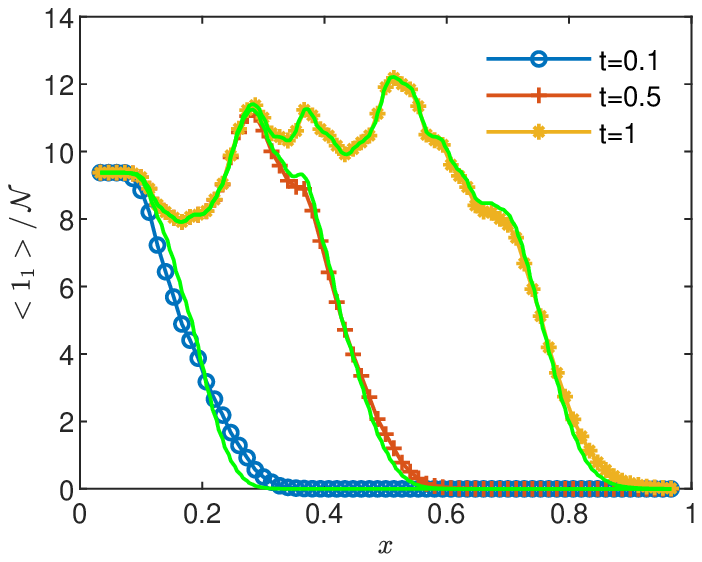}
\caption{\label{fig:rand_CGc1}Coarse-grained contaminant concentrations compared to volume averages (full lines).}
\end{minipage}
\hspace*{0.1in}
\begin{minipage}[t]{0.45\linewidth}\centering
\includegraphics[width=\linewidth]{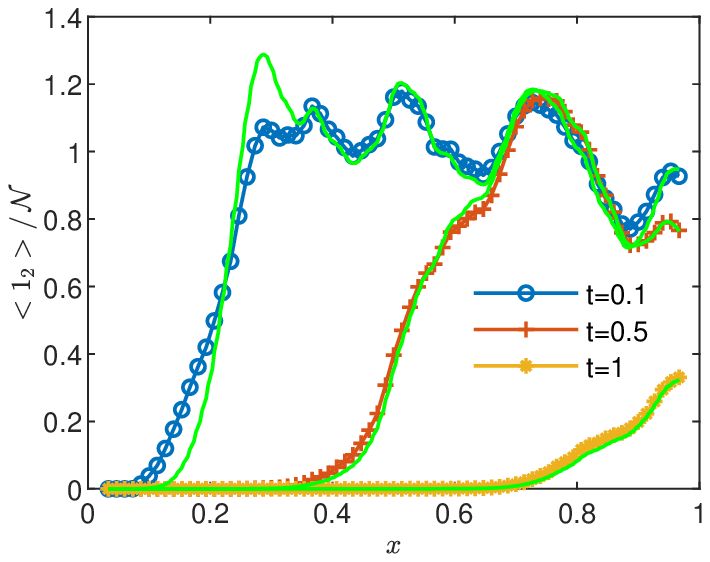}
\caption{\label{fig:rand_CGc2}Coarse-grained oxygen concentrations compared to volume averages (full lines).}
\end{minipage}
\end{figure}

\begin{figure}
\begin{minipage}[t]{0.45\linewidth}\centering
\includegraphics[width=\linewidth]{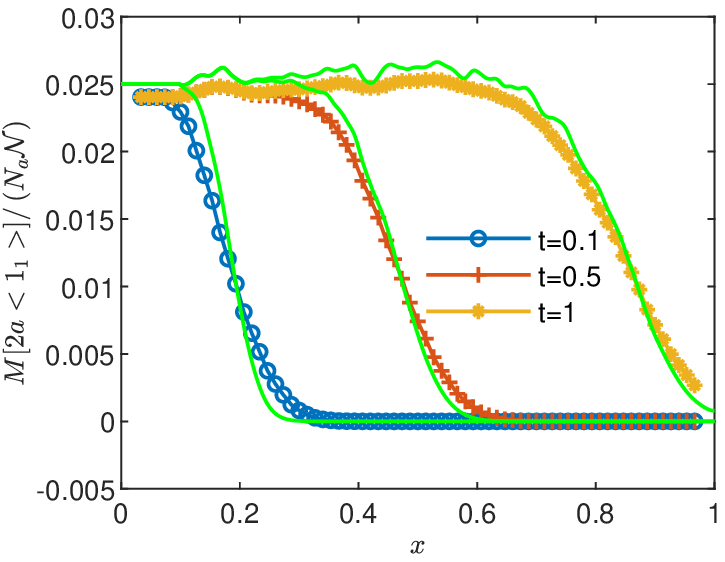}
\caption{\label{fig:MC_c1}Ensemble averaged moving averages and fine-grained contaminant concentrations $c_{1}$ (full lines).}
\end{minipage}
\hspace*{0.1in}
\begin{minipage}[t]{0.45\linewidth}\centering
\includegraphics[width=\linewidth]{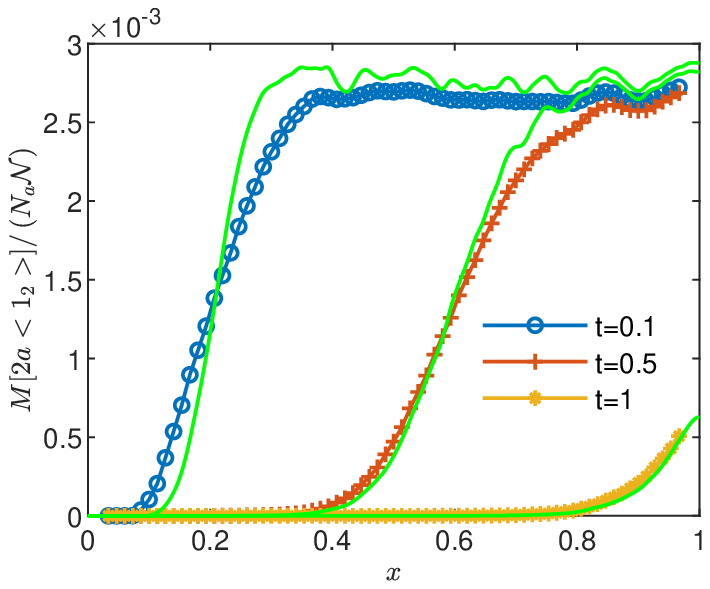}
\caption{\label{fig:MC_c2}Ensemble averaged moving averages and fine-grained oxygen concentrations $c_{2}$ (full lines).}
\end{minipage}
\end{figure}

\begin{figure}
\begin{minipage}[t]{0.45\linewidth}\centering
\includegraphics[width=\linewidth]{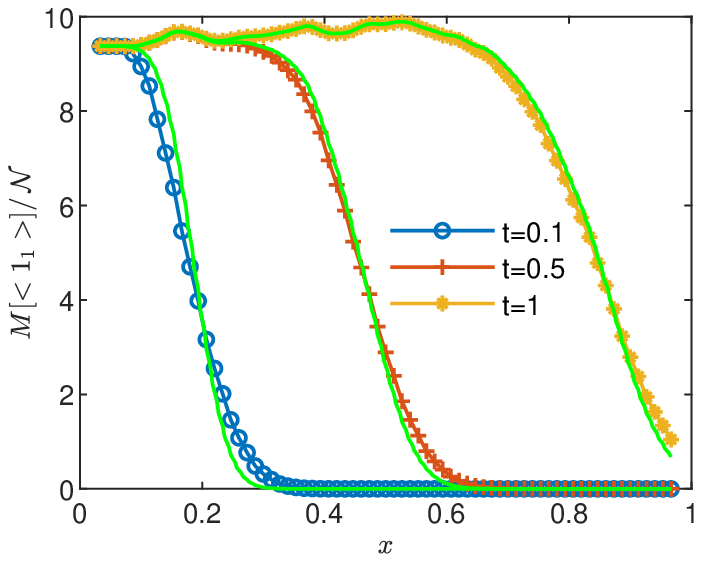}
\caption{\label{fig:MC_CGc1}Ensemble averaged coarse-grained and volume averages (full lines) of the contaminant concentrations.}
\end{minipage}
\hspace*{0.1in}
\begin{minipage}[t]{0.45\linewidth}\centering
\includegraphics[width=\linewidth]{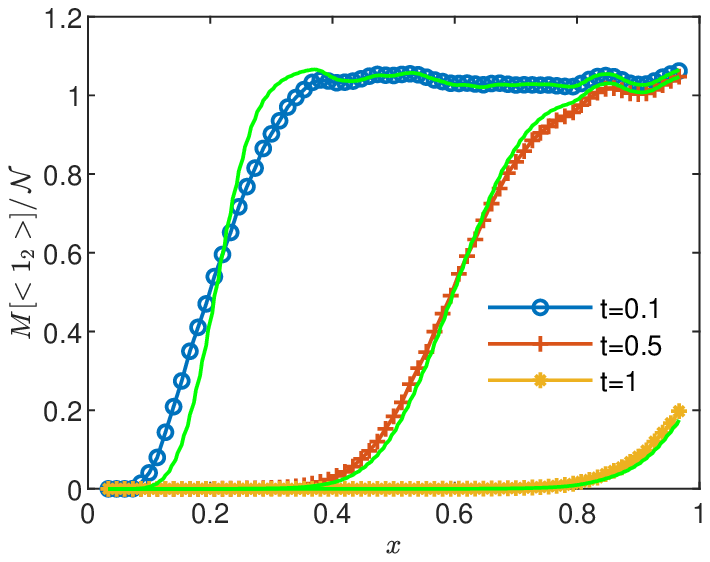}
\caption{\label{fig:MC_CGc2}Ensemble averaged coarse-grained and volume averages (full lines) of the oxygen concentrations.}
\end{minipage}
\end{figure}

The initial contaminant concentration is set to one mole uniformly distributed on the lattice sites $l=1,\ldots,\Delta \mathcal{L}$, with $\Delta \mathcal{L}=(\mathcal{L}-1)/10$, $\mathcal{L}=1/\Delta x +1$, and zero otherwise. Complementary, the oxygen concentration is set to one mole uniformly distributed on the lattice sites $l=\Delta \mathcal{L}+1,\ldots, \mathcal{L}$ and zero otherwise. The condition $n_{\mathcal{L},k}=n_{\mathcal{L}-1,k}$ is imposed on the outflow boundary and the first $\Delta \mathcal{L}$ sites are set to the initial condition after each time step $k$. One simulates in this way a constant injection of contaminant near the inflow boundary, the consumption of the contaminant and oxygen through the biodegradation process in $\Omega$, and the evacuation of the mixture of contaminant and oxygen at the outflow boundary. The solution is obtained with the same scheme as in the previous section. The final time is set to $T=1.1$ and the spatial and temporal scales are $a=0.03$ and $\tau=0.1$. The discretization parameters consistent with the BGRW algorithm are $\Delta x=2.50\cdot 10^{-3}$ and $\Delta t=1.56\cdot 10^{-4}$. In Figs.~\ref{fig:rand_c1}~-~\ref{fig:MC_CGc2} we present the results obtained for both single-realizations and for an ensemble of 100 realizations of the random velocity field.

Figures~\ref{fig:rand_c1}~and~\ref{fig:rand_c2} show the moving averages and the actual concentration of the two species at the selected sampling times, for the same realization of the velocity field. The coarse-grained averages are compared to the volume averages in Figs.~\ref{fig:rand_CGc1}~and~\ref{fig:rand_CGc2}. One remarks that the two averages can be significantly different mainly at the first sampling time. We also remark the much larger magnitude of the space-time upscaled concentrations as compared to the fine-grained concentrations from Figs.~\ref{fig:rand_c1}~and~\ref{fig:rand_c2}. In Figs.~\ref{fig:MC_c1}~and~\ref{fig:MC_c2} we can see that the ensemble averaging smooths out the moving averages but not the actual concentrations, which still show oscillations after averaging over 100 realizations. The averaged concentrations remain in the same range as those for single realizations. Thus, the ensemble averaging, even combined with moving averages, does not perform an upscaling of the fine-grained concentration. Instead, as seen in Figs.~\ref{fig:MC_CGc1}~and~\ref{fig:MC_CGc2}, the ensemble averaging smooths out the upscaled concentration obtained by CGST averaging, the result being a stochastic upscaling adapted to the spatial scale $a$ and the temporal scale $\tau$ on an hypothetic measurement. Similarly to the single realization results, differences between the ensemble averaged CGST and volume averages are visible at the first sampling time. Instead, the two averages are close to each other in regions without reactions between contaminant and oxygen, e.g., contaminant concentration for $x<0.5$ at $t=1$ (Figs.~\ref{fig:MC_CGc1}~and~\ref{fig:rand_CGc1}) and oxygen concentration for $x>0.5$ at $t=0.1$ (Figs.~\ref{fig:MC_CGc2}~and~\ref{fig:rand_CGc2}). A more precise evaluation of the discrepancy is obtained by computing relative differences. Table~\ref{table:MaxErr_1Daquifer} shows that while the global differences $e_{c_{\nu}}$, computed with $l^2$-norms are generally smaller than 10\%, for both single realization and ensemble averages, the maximum relative differences $\varepsilon_{c_{\nu}}$ are larger than 13\% at all the three sampling times, up to 53\% for ensemble averages, at $t=0.1$. Thus, even though the flow is stationary, time-dependent reactions make the transport process nonstationary and the volume average alone could be an inaccurate model of the measurement.

\begin{table}
\caption{\label{table:MaxErr_1Daquifer}Relative differences between volume and CGST averages: biodegradation in one-dimensional aquifer with random hydraulic conductivity (single realization and ensemble averages).}
\begin{center}
\begin{tabular}{ c | c c| c c | c c| c c c }
 & \multicolumn{4}{c|}{realization} & \multicolumn{4}{c}{ensemble} \\
\hline
    $t$ & $e_{c_1}$ & $\varepsilon_{c_1}$ & $e_{c_2}$ & $\varepsilon_{c_2}$ & $e_{c_1}$ & $\varepsilon_{c_1}$ & $e_{c_2}$ & $\varepsilon_{c_2}$ & \\
  \hline
     0.1  & 0.1018 & 0.1888 & 0.0848 & 0.2110 & 0.0969 & 0.1718 & 0.0607 & 0.5275 \\
     0.5  & 0.0277 & 0.3186 & 0.0323 & 0.1605 & 0.0284 & 0.3098 & 0.0347 & 0.1809 \\
     1.0  & 0.0207 & 0.3143 & 0.0510 & 0.1515 & 0.0161 & 0.2768 & 0.1775 & 0.1365 \\
  \hline
\end{tabular}
\end{center}
\end{table}

\begin{remark}
\label{rem:MCaverage}
In simulations of transport in random velocity fields, a quantity of practical importance is the ensemble averaged concentration. The stochastic average $M[\cdot]$ used in Eq.~(\ref{eq:stoch_average}) also estimates this quantity if the mean is computed over all the realizations of the random walk process, as well as over all the realizations of the random velocity. The GRW algorithms are self-averaging, in the sense that if the total number of particles $\mathcal{N}$ is large enough no averages over the realizations of the random walk process are needed to obtain the desired accuracy. For instance, in \cite[Fig. 5]{Vamosetal2003} it is shown that the GRW solution for a one-dimensional diffusion problem similar to those considered in this study reaches the accuracy of the equivalent finite difference scheme in a single realization of the random walk by using $10^{6}$ particles. Since in the current simulations $\mathcal{N}$ is of the order of Avogadro's number, it follows that the averages $M[\langle 1_{\nu} \rangle]$ presented in Figs.~\ref{fig:MC_CGc1}~and~\ref{fig:MC_CGc2} correspond to the space-time averages (\ref{eq:stoch_average}) of the ensemble averaged fine-grained concentrations $c_{\nu}$ plotted with full lines in Figs.~\ref{fig:MC_c1}~and~\ref{fig:MC_c2}.
\end{remark}

\subsubsection{Monod reactions in variably saturated soils}
\label{sec:MonodReaction_soil}

One expects that for time-dependent reactions in non stationary flow regimes the discrepancy between volume and CGST averages will be even more pronounced. As an example, we consider a one-dimensional model problem for biodegradation processes in soil columns. The water flow in variably saturated porous media is described by the one-dimensional Richards equation
\begin{equation}\label{eq:Richards}
\partial_t\theta(\psi)-\partial_z\left[K(\theta(\psi))\partial_z(\psi+z)\right]=0,
\end{equation}
where $\psi(z,t)$ is the pressure head expressed in length units, $\theta$ is the volumetric water content, $K$ stands for the hydraulic conductivity of the medium, and $z$ is the height oriented positively upward. The water flux given by Darcy's law is $q=-K(\theta(\psi))\frac{\partial}{\partial z}(\psi+z)$. The reactive transport is now governed by Richards equation (\ref{eq:Richards}) coupled with the transport equations
\begin{equation}\label{eq:reactRichards}
\partial_t\left[\theta(\psi) c_{\nu}\right] + \partial_z(qc_{\nu})-D\partial^{2}_{z}c_{\nu}=R_{\nu}(c_{1},c_{2},\theta),\; \nu=1,2.
\end{equation}

We consider in the following the relationships defining the water content $\theta(\psi)$ and the hydraulic conductivity $K(\theta(\psi))$ given by the van Genuchten-Mualem model
\begin{equation} \label{eq:theta}
\Theta(\psi) = \begin{cases} \left(1+(-\alpha \psi)^n\right)^{-m}, &\psi < 0 \\
1, &\psi \geq 0,
\end{cases}
\end{equation}
\begin{equation} \label{eq:K}
K(\Theta(\psi)) = \begin{cases} K_{sat} \Theta(\psi)^{\frac{1}{2}} \left[1-\left(1-\Theta(\psi)^\frac{1}{m}\right)^m \right]^2, &\psi < 0 \\
K_{sat}, &\psi \geq 0,
\end{cases}
\end{equation}
where $\Theta = (\theta - \theta_{res})/(\theta_{sat} - \theta_{res})$ is the normalized water content, $\theta_{res}$ and $\theta_{sat}$ are the residual and the saturated water content, $K_{sat}$ is the hydraulic conductivity of the saturated soil, and $\alpha$, $n$, $m=1-1/n$ are model parameters depending on the soil type. In this example, we consider the parameters of a silt loam used in previously published studies \cite{ListandRadu2016,Illiano2020,Suciuetal2021,SuciuandRadu2021}, $\theta_{sat}=0.396$, $\theta_{res}=0.131$, $\alpha=0.423$, $n=2.06$, $K_{sat}=4.96 \cdot 10^{-2}$. The soil heterogeneity is modeled by a fixed realization of a random saturated hydraulic conductivity with mean equal to the parameter $K_{sat}$ of the loam soil model. The realization is computed with the Kraichnan method described in \cite[Appendix C.3.1.2]{Suciu2019}, for a random $\ln K$ field with Gaussian correlation specified by the correlation length $\lambda=0.1$ and the variance $\sigma^2=0.5$. The dispersion is parameterized by the constant coefficient $D=0.001$. The system of coupled equations (\ref{eq:Richards}-\ref{eq:reactRichards}), parameterized by the model (\ref{eq:theta}-\ref{eq:K}), with reactions governed by the Monod model (\ref{eq:monod}) using the same parameters as in the saturated case presented above, is solved in the domain $\Omega=(0,3)$ for the final time $T=6.6$. The space-time scale of the CGST procedure is specified by $a=0.1$ and $\tau=0.6$.

The infiltration of the contaminated water in the column is driven by the variable pressure on the top boundary given by $\psi(3,t\leq t_1)=-3+3.2 t/t_1$ and $\psi(3,t> t_1)=0.2$, with $t_1=T/3$. The initial pressure corresponds to the unsaturated profile $\psi(z,0)=-z$ and the constant pressure condition $\psi(0,t)=\psi(0,0)$ is set on the bottom boundary. The initial contaminant concentration is set to one mole uniformly distributed on the lattice sites $l=\mathcal{L}-\Delta \mathcal{L},\ldots, \mathcal{L}$, with $\Delta \mathcal{L}=(\mathcal{L}-1)/10$, $\mathcal{L}=3/\Delta z +1$, on the top of the domain $\Omega$ and zero otherwise. The oxygen concentration is set to one mole uniformly distributed on the lattice sites $l=1,\ldots, \mathcal{L}-\Delta \mathcal{L}-1$ and zero otherwise. The top $\Delta \mathcal{L}$ sites are set to the initial condition after each time step $k$ and the condition $n_{1,k}=n_{2,k}$ is imposed on the bottom boundary for both species. With these, one simulates a continuous injection of contaminant over a length $\Delta \mathcal{L}\Delta z$ on the top of the column, assuming that the incoming contaminated water does not contain oxygen, in conditions of free drainage at the bottom of the column.

The flow equation (\ref{eq:Richards}) is solved with the iterative GRW algorithm for one-dimensional flows, implemented as an $L$-scheme linearization approach (see \cite{Popetal2004,ListandRadu2016}) capable to handle the nonlinearity and the degeneracy of the Richards equation \cite[Sect. 2]{Suciuetal2021}. At each time step $k$, the solution $(\psi, \theta)$ of the Richards equation is transmitted to the reactive transport solver. The
latter is an iterative $L$-scheme, which at each iteration step solves the system (\ref{eq:reactRichards}) with an operator splitting procedure consisting of a BGRW solution of the advection-diffusion process and a deterministic solution of the reaction step (this is the one-dimensional version of the two-dimensional iterative scheme introduced in \cite[Sect. 3.2]{SuciuandRadu2021}). The spatial step is fixed to $\Delta z=1.25\cdot 10^{-2}$. The time steps, determined as the minimum of the time steps required by the GRW flow solver and the BGRW transport solver lay between $\Delta t=7.42\cdot 10^{-3}$ and $\Delta t=7.3\cdot 10^{-3}$. The number of iterations required for the convergence of the two $L$-schemes ranges between 100 and 450.

\begin{figure}
\begin{minipage}[t]{0.45\linewidth}\centering
\includegraphics[width=\linewidth]{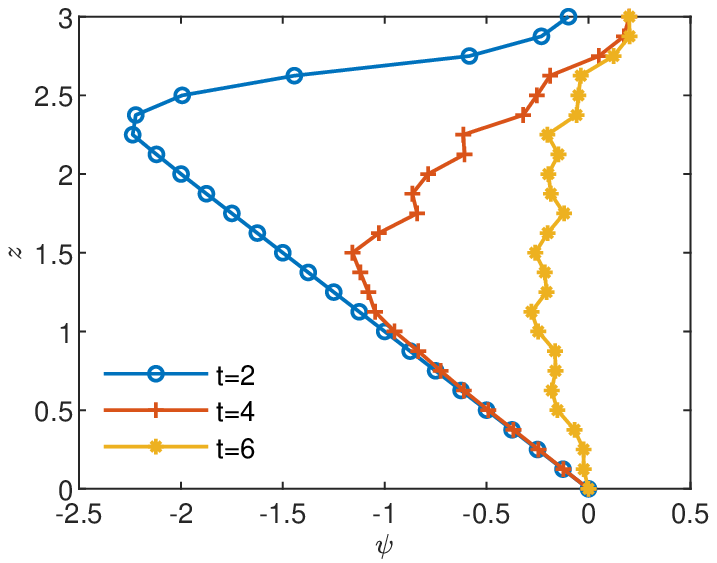}
\caption{\label{fig:p}Pressure profiles through the soil column at three sampling times.}
\end{minipage}
\hspace*{0.1in}
\begin{minipage}[t]{0.45\linewidth}\centering
\includegraphics[width=\linewidth]{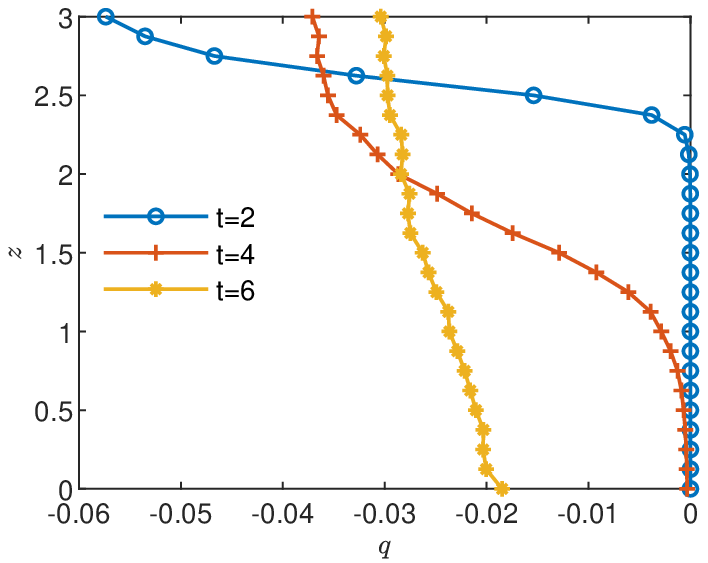}
\caption{\label{fig:q}Flow velocity profiles through the soil column at three sampling times.}
\end{minipage}
\end{figure}

\begin{figure}
\begin{minipage}[t]{0.45\linewidth}\centering
\includegraphics[width=\linewidth]{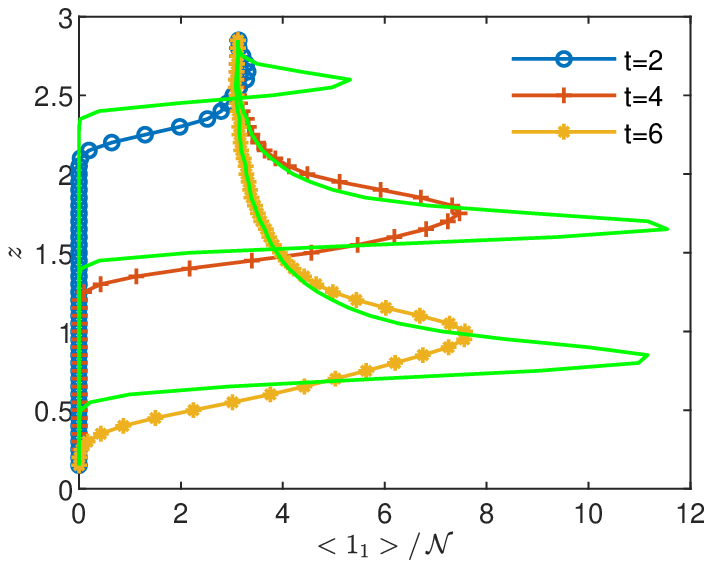}
\caption{\label{fig:CGc1}Profiles of the CGST upscaled contaminant concentration at three sampling times compared to volume averages (full lines).}
\end{minipage}
\hspace*{0.1in}
\begin{minipage}[t]{0.45\linewidth}\centering
\includegraphics[width=\linewidth]{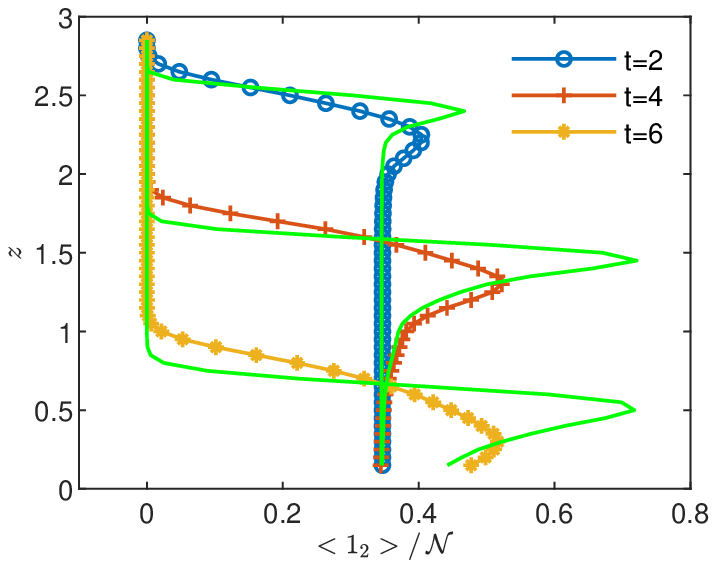}
\caption{\label{fig:CGc2}Profiles of the CGST upscaled oxygen concentration at three sampling times  compared to volume averages (full lines)}
\end{minipage}
\end{figure}

\begin{table}
\caption{\label{table:MaxErr_1Dsoil}Relative differences between volume and CGST averages: one-dimensional model of soil column, random saturated hydraulic conductivity, unsaturated/saturated flow regimes.}
\begin{center}
\begin{tabular}{ c | c c | c c }
   $ t $ & $e_{c_1}$ & $\varepsilon_{c_1}$ &  $e_{c_2}$ & $\varepsilon_{c_2}$ \\
\hline
     2 & 0.5054 & 0.9895 & 0.1127 & 0.5867 \\
     4 & 0.3420 & 0.6950 & 0.2392 & 0.6057 \\
     6 & 0.3136 & 0.7702 & 0.3745 & 0.6576 \\
  \hline
\end{tabular}
\end{center}
\end{table}

The flow solutions and the space-time upscaled concentrations are presented in Figs.~\ref{fig:p}-\ref{fig:CGc2}. The pressure profiles recorded at the three sampling times shown in Fig.~\ref{fig:p} capture the spatial heterogeneity of the flow solution and the transition from the unsaturated regime ($\psi<0$) to the saturated regime ($\psi\geq 0$). The resulting flow velocity (Fig.~\ref{fig:q}) is strongly variable in both space and time. The complexity of the unsaturated/saturated flow impacts on the distribution of the reacting species (see Figs.~\ref{fig:CGc1}~and~\ref{fig:CGc2}), which inherits the heterogeneity and the spatio-temporal variability of the flow. In these conditions, the discrepancy between the volume and the CGST average is dramatically large at all the tree sampling times, with maximum relative differences between 59\% and 99\% (Table~\ref{table:MaxErr_1Dsoil}). Thus, averaging over the spatial interval $(x-a,x+a)$ at given time $t$ is completely irrelevant for the observation/measurement made at the spatio-temporal scale $(a,\tau)$ centered at $(x,t)$. Excepting special situations, such as the case of almost homogeneous reactions with slow time variations presented in Figs.~\ref{fig:bimol_CGc1}~and~\ref{fig:bimol_CGc2}, the time average has to be considered as well in one-dimensional modeling of experimental observations.

\section{CGST averaging in two-dimensional transport problems}
\label{sec:cgst_2dim}

The two-dimensional model for reactive transport with biodegradation reactions, similar to the one-dimensional model considered above, is governed by the equations
\begin{align}
&\partial_t\theta(\psi)-\nabla\cdot\left[K(\theta(\psi)\nabla(\psi+z)\right]=0,\label{eq:Richards2D}\\
&\partial_t\left[\theta(\psi) c_{\nu}\right]-\nabla\cdot\left(D\nabla c_{\nu} -\mathbf{q}c_{\nu}\right)=R_{\nu}(c_{1},c_{2},\theta),\;\;\mathbf{q}=-K(\theta(\psi)\nabla(\psi+z), \;\; \nu=1,2.\label{eq:reactRichards2D}
\end{align}
Equations~(\ref{eq:Richards2D}-\ref{eq:reactRichards2D}) are parameterized with the van Genuchten-Mualem model (\ref{eq:theta}-\ref{eq:K}) and the same Monod reaction system as in the previous sections. The two-dimensional CGST averaging procedure, first verified by computing the velocity components and the diffusion coefficients, will be used in the following to compute space-time upscaled concentrations.

\subsection{Verification of the two-dimensional CGST averaging procedure}

The two-dimensional CGST averaging procedure corresponds to $d=2$ in the general definition (\ref{eq:average}). Similarly to the derivation of the one-dimensional stationary diffusion equation (\ref{eq:diff}) given in Section~\ref{sec:verification}, one obtains the diffusion equation
\begin{equation*}
\partial_{x_{\alpha}}\partial_{x_{\beta}}(D_{{\alpha\beta}}\;c)=0,
\end{equation*}
with diffusion tensor defined by $D_{{\alpha\beta}}=-\left\langle x_{\alpha}\xi_{\beta}\right\rangle/\left\langle 1\right\rangle$. The numerical implementation of the CGST averaging procedure is a straightforward extension of that used in the one-dimensional case.

We perform BGRW/GRW simulations in the two-dimensional domain $\Omega=(0,1)\times(0,1)$ of the advection-diffusion process defined by constant diagonal diffusion tensor with components $D_{11}=D_{22}=0.0001$ and constant velocity $(u,v)=(1,0)$. The BGRW/GRW simulations are done with the same number of particles and space-time discretization and the CGST averages are computed with the same parameters $a=0.05$ and $\tau=0.1$ as in the one dimensional case. The components of the diffusion tensor (Tables~\ref{table:Dbgrw}~and~\ref{table:Dgrw}) and of the velocity (Table~\ref{table:uv}) computed with the CGST procedure are, again, in very good agreement with the nominal coefficients of the advection-diffusion process.

\begin{table}
\caption{\label{table:Dbgrw}CGST computation of the intrinsic diffusion coefficients: BGRW algorithm.}
\begin{center}
\begin{tabular}{ c c c c c }
    $\Delta x$ & $\overline{D_{11}}$ & $\overline{D_{12}}$  & $\overline{D_{21}}$ & $\overline{D_{22}}$ \\
  \hline
     5.00e-03  & 0.0001$\pm$ 1.59e-19 & -1.59e-20$\pm$ 8.13e-20 & 0.00$\pm$ 0.00 & \hspace{-0.6cm}0.0001$\pm$ 0.00 \\
     2.50e-03  & 0.0001$\pm$ 1.31e-19 & -1.06e-19$\pm$ 2.92e-19 & 0.00$\pm$ 0.00 & \hspace{-0.6cm}0.0001$\pm$ 0.00 \\
     1.25e-03  & 0.0001$\pm$ 3.02e-19 & \hspace{0.15cm}5.93e-20$\pm$ 2.36e-19 & 0.00$\pm$ 0.00 & 0.0001$\pm$ 2.87e-20 \\
     6.25e-04  & 0.0001$\pm$ 7.47e-19 & \hspace{0.15cm}1.23e-19$\pm$ 3.02e-19 & 0.00$\pm$ 0.00 & 0.0001$\pm$ 1.44e-20 \\
  \hline
\end{tabular}
\end{center}
\end{table}

\begin{table}
\caption{\label{table:Dgrw}CGST computation of the intrinsic diffusion coefficients: GRW algorithm.}
\begin{center}
\begin{tabular}{ c c c c c }
    $\Delta x$ & $\overline{D_{11}}$ & $\overline{D_{12}}$  & $\overline{D_{21}}$ & $\overline{D_{22}}$ \\
  \hline
     5.00e-03  & 0.0001$\pm$ 1.20e-19 & -3.57e-21$\pm$ 7.29e-20 & 0.00$\pm$ 0.00 & 0.0001$\pm$ 1.44e-20 \\
     2.50e-03  & 0.0001$\pm$ 1.63e-19 & -5.21e-20$\pm$ 2.63e-19 & 0.00$\pm$ 0.00 & 0.0001$\pm$ 2.87e-20 \\
     1.25e-03  & 0.0001$\pm$ 2.46e-19 & -1.09e-19$\pm$ 2.73e-19 & 0.00$\pm$ 0.00 & 0.0001$\pm$ 1.44e-20 \\
     6.25e-04  & 0.0001$\pm$ 9.66e-19 & \hspace{0.15cm}2.07e-19$\pm$ 3.56e-19 & 0.00$\pm$ 0.00 & 0.0001$\pm$ 1.44e-20 \\
  \hline
\end{tabular}
\end{center}
\end{table}

\begin{table}
\caption{\label{table:uv}CGST computation of the velocity components.}
\begin{center}
\begin{tabular}{ c c c | c c c }
    & \hspace{2cm}$\overline{u}$  & & \hspace{2cm}$\overline{v}$  \\
\hline
    $\Delta x$ & BGRW & GRW  & BGRW & GRW \\
  \hline
     5.00e-03  & \hspace{-0.6cm}1.00$\pm$ 0.00 & 1.00$\pm$ 0.00 & \hspace{-0.6cm}-2.71e-19$\pm$ 0.00 & 0.00$\pm$ 0.00 \\
     2.50e-03  & 1.00$\pm$ 1.18e-16 & 1.00$\pm$ 0.00 & \hspace{-0.6cm}-1.39e-19$\pm$ 0.00 & 0.00$\pm$ 0.00 \\
     1.25e-03  & \hspace{-0.6cm}1.00$\pm$ 0.00 & 1.00$\pm$ 0.00 & \hspace{0.15cm}1.14e-19$\pm$ 2.55e-35 & 0.00$\pm$ 0.00 \\
     6.25e-04  & 1.00$\pm$ 2.36e-16 & 1.00$\pm$ 0.00 & -9.55e-20$\pm$ 1.28e-35 & 0.00$\pm$ 0.00 \\
  \hline
\end{tabular}
\end{center}
\end{table}

\subsection{Application to biodegradation in variably saturated soils}
\label{sec:MonodSoil2dim}

\begin{figure}
\includegraphics[width=\linewidth]{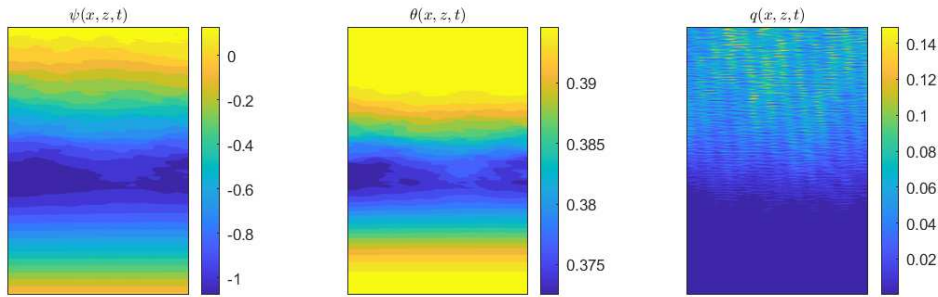}
\caption{\label{fig:T120_p}Pressure $\psi$, water content $\theta$, and velocity amplitude $q=\sqrt{q_{x}^2+q_{z}^2}$ at the final time $T=132$.}
\end{figure}

\begin{figure}
\begin{minipage}[t]{0.49\linewidth}\centering
\includegraphics[width=\linewidth]{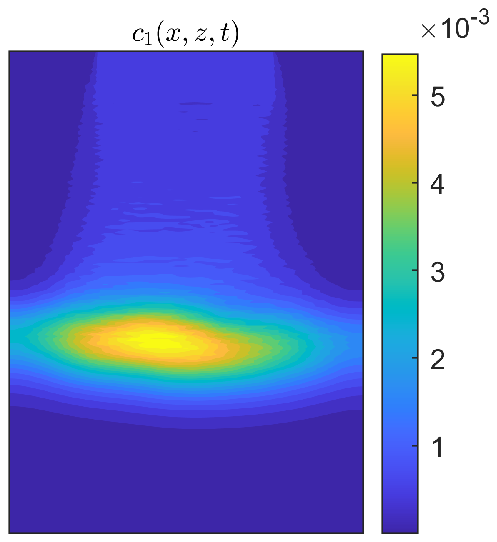}
\caption{\label{fig:T120_c1}Contaminant concentration $c_1$ at the final time $T=132$.}
\end{minipage}
\hspace*{0.1in}
\begin{minipage}[t]{0.49\linewidth}\centering
\includegraphics[width=\linewidth]{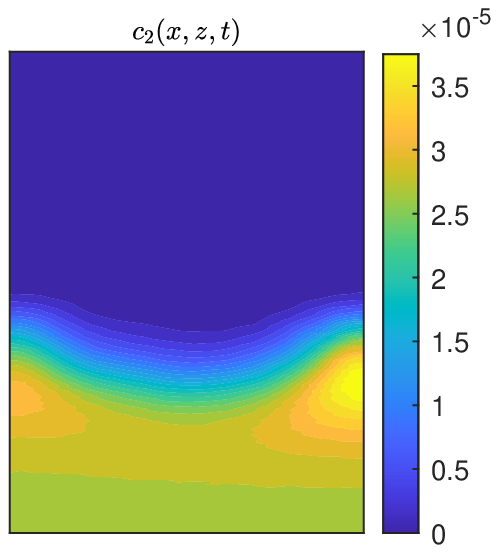}
\caption{\label{fig:T120_c2}Contaminant concentration $c_1$ at the final time $T=132$.}
\end{minipage}
\end{figure}

\begin{figure}
\begin{minipage}[t]{0.45\linewidth}\centering
\includegraphics[width=\linewidth]{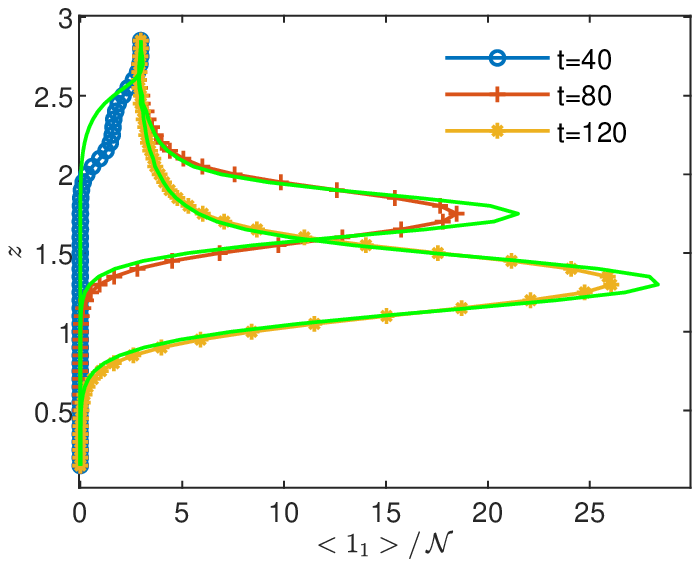}
\caption{\label{fig:T120_cg1_centr}Profiles of the CGST upscaled contaminant concentration at three sampling times compared to volume averages (full lines), sampled on a centered vertical line.}
\end{minipage}
\hspace*{0.1in}
\begin{minipage}[t]{0.45\linewidth}\centering
\includegraphics[width=\linewidth]{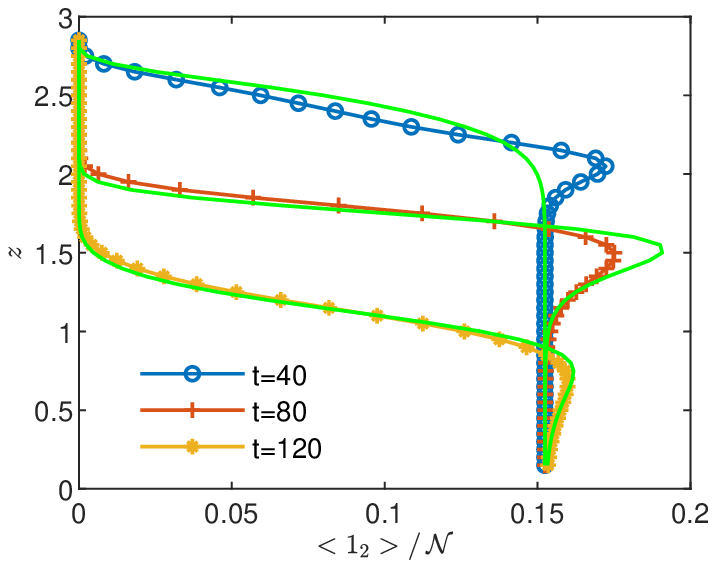}
\caption{\label{fig:T120_cg2_centr}Profiles of the CGST upscaled oxygen concentration at three sampling times compared to volume averages (full lines), sampled on a centered vertical line through the domain.}
\end{minipage}
\end{figure}

To illustrate the two-dimensional CGST averaging procedure, we simulate the same biodegradation process as in the one-dimensional case in the domain $\Omega=(0,2)\times(0,3)$. We consider the same soil model, vertical profile of the pressure, the same top and bottom pressure boundary conditions, and no-flow conditions on the left and right boundaries. The initial concentrations for contaminant and oxygen in the rectangle $\Omega_1=[0.5,1.5]\times[2.7,3]$ and in the surrounding region take the same values as on the top $\Delta\mathcal{L}$ sites and on the remaining sites of the lattice in the one-dimensional problem and concentrations in $\Omega_1$ are set to the initial condition after each time step. No-flux conditions are imposed on the bottom and left/right side boundaries for both the contaminant and the oxygen concentrations. CGST averages are computed with $a=0.1$ and $\tau=12$, for a total simulation time $T=132$. The flow and the BGRW solutions are computed to the flow $L$-scheme from \cite[Sect. 4.1]{Suciuetal2021} coupled with the iterative splitting procedure introduced in \cite[Sect. 3.2]{SuciuandRadu2021}.

\begin{figure}
\begin{minipage}[t]{0.45\linewidth}\centering
\includegraphics[width=\linewidth]{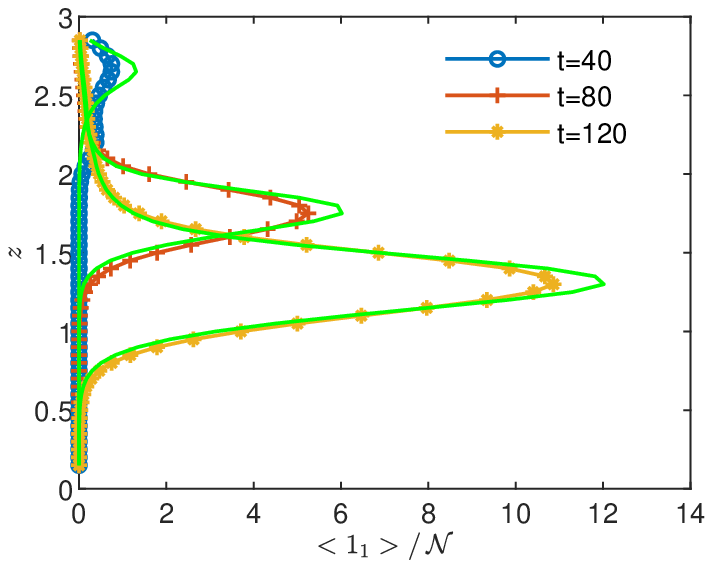}
\caption{\label{fig:T120_cg1_dec}Profiles of the CGST upscaled contaminant concentration at three sampling times  compared to volume averages (full lines), sampled on a vertical line located at $x=1.75$}
\end{minipage}
\hspace*{0.1in}
\begin{minipage}[t]{0.45\linewidth}\centering
\includegraphics[width=\linewidth]{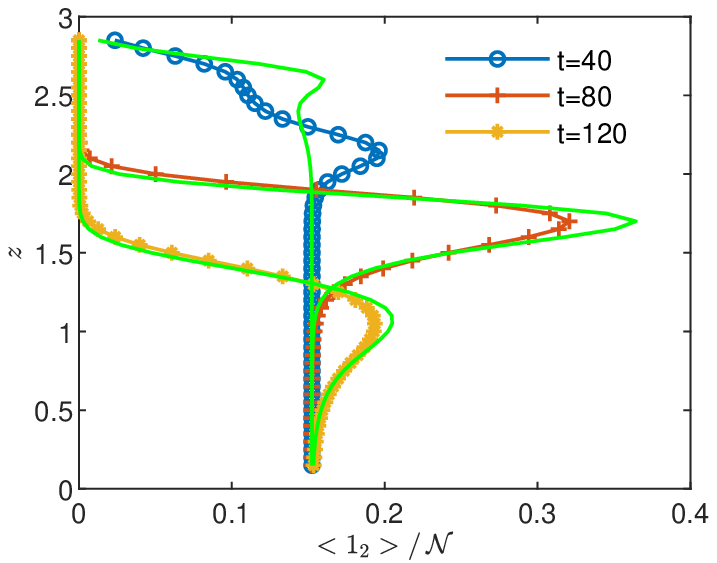}
\caption{\label{fig:T120_cg2_dec}Profiles of the CGST upscaled oxygen concentration at three sampling times  compared to volume averages (full lines), sampled on a vertical line located at $x=1.75$}
\end{minipage}
\end{figure}

CGST averages are computed along a vertical sampling line placed at the center of the horizontal edge, as well as along a vertical line close to the right boundary of the domain. Figures~\ref{fig:T120_p}~-~\ref{fig:T120_c2} show the flow solution and the fine-grained contaminant and oxygen concentrations at the final time $T=132$. Significant differences between the CGST and volume averages at the three sampling times can be observed for both contaminant and oxygen concentrations presented in Figs.~\ref{fig:T120_cg1_centr}-\ref{fig:T120_cg2_dec} in regions where the two species interact and negligible differences outside the reaction zone (indicated by constant concentrations). One remarks also that in case of decentered sampling, near the right boundary (Figs.~\ref{fig:T120_cg1_dec}~and~\ref{fig:T120_cg2_dec}), these differences are generally larger than for centered sampling line (Figs.~\ref{fig:T120_cg1_centr}~and~\ref{fig:T120_cg2_centr}). Table~\ref{table:MaxErr_2Dsoil} quantifies these visual observations and warns that disregarding the time average when modeling measurements through two-dimensional simulations of reactive transport in partially saturated soils may result in discrepancies up to 80\% with respect to the CGST averaging approach.

\begin{table}
\caption{\label{table:MaxErr_2Dsoil}Relative differences between volume and CGST averages: two-dimensional model of soil column, random saturated hydraulic conductivity, unsaturated/saturated flow regimes, centered ($x=1$) and decentered ($x=1.75$) vertical sampling line.}
\begin{center}
\begin{tabular}{ c | c c | c c | c c | c c c }
 & \multicolumn{4}{c|}{centered} & \multicolumn{4}{c}{decentered} \\
\hline
    $ t $ & $e_{c_1}$ & $\varepsilon_{c_1}$ & $e_{c_2}$ & $\varepsilon_{c_2}$ & $e_{c_1}$  & $\varepsilon_{c_1}$ & $e_{c_2}$ & $\varepsilon_{c_2}$ \\
  \hline
     40  & 0.3466 & 0.8041 & 0.0675 & 0.1389 & 0.6080 & 0.7670 & 0.1320 & 0.5471 \\
     80  & 0.1198 & 0.1636 & 0.0536 & 0.3479 & 0.1426 & 0.1737 & 0.0776 & 0.1354 \\
     120 & 0.0654 & 0.0880 & 0.0280 & 0.0462 & 0.0862 & 0.1062 & 0.0440 & 0.2204 \\
  \hline
\end{tabular}
\end{center}
\end{table}

\subsection{Application to biodegradation in saturated aquifers}

To complete the picture of the CGST averaging in two dimensions, we present the simulation of the biodegradation process in conditions of a continuous injection at the inflow boundary of a saturated aquifers contained in the domain $\Omega=(0,20)\times(0,10)$ and evacuation of of the solute at the outflow boundary, similar to that considered in \cite{SuciuandRadu2021}. The two-dimensional velocity field is approximated with the Kraichnan procedure described in Section~\ref{sec:MonodReaction_randV} above, for a $\ln(K)$ field with correlation length $\lambda=1$ m and variance $\sigma^2=0.5$. The 2.5 moles of contaminant are continuously injected in a square $0.5\times0.5$ on the middle of the inflow boundary and the initial oxygen concentration is set to 0.1 moles outside the injection domain. No-flux conditions are imposed on the outflow boundary while the concentrations on left, top, and bottom boundaries are kept at the initial values. The biodegradation is governed by the same Monod model as in the previous sections. CGST averages are computed with $a=0.67$ and $\tau=30$, for a total simulation time $T=330$. The BGRW solution is computed with the iterative splitting procedure from \cite[Sect. 3.2]{SuciuandRadu2021}.

\begin{figure}
\begin{minipage}[t]{0.49\linewidth}\centering
\includegraphics[width=\linewidth]{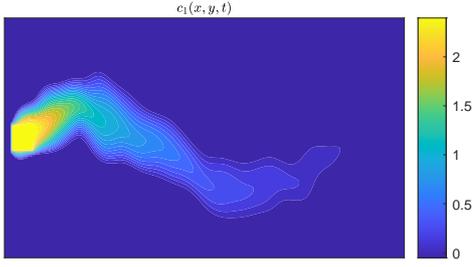}
\caption{\label{fig:aqv_c1}Contaminant concentration at  $T=330$.}
\end{minipage}
\hspace*{0.1in}
\begin{minipage}[t]{0.49\linewidth}\centering
\includegraphics[width=\linewidth]{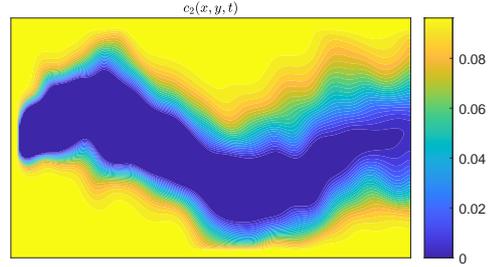}
\caption{\label{fig:aqv_c2}Contaminant concentration at $T=330$.}
\end{minipage}
\end{figure}

\begin{figure}
\begin{minipage}[t]{0.45\linewidth}\centering
\includegraphics[width=\linewidth]{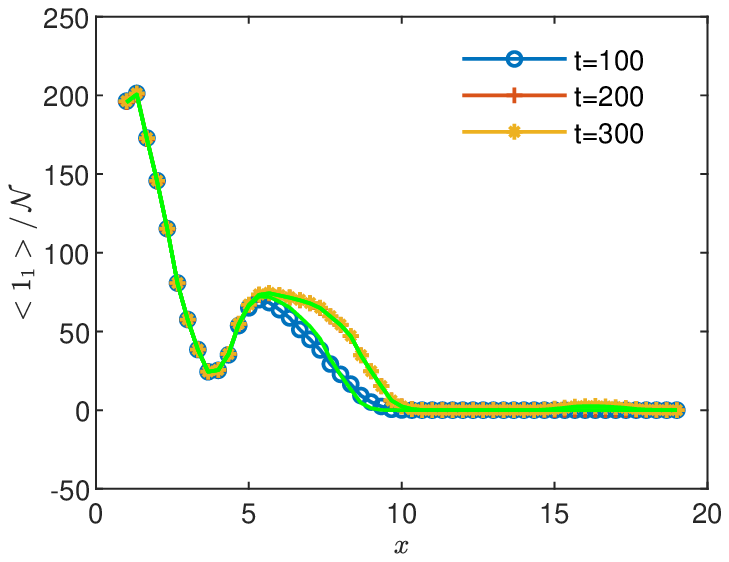}
\caption{\label{fig:aqv_cg1}Profiles of the CGST upscaled contaminant concentration at three sampling times compared to volume averages (full lines), sampled on the horizontal midline of the domain.}
\end{minipage}
\hspace*{0.1in}
\begin{minipage}[t]{0.45\linewidth}\centering
\includegraphics[width=\linewidth]{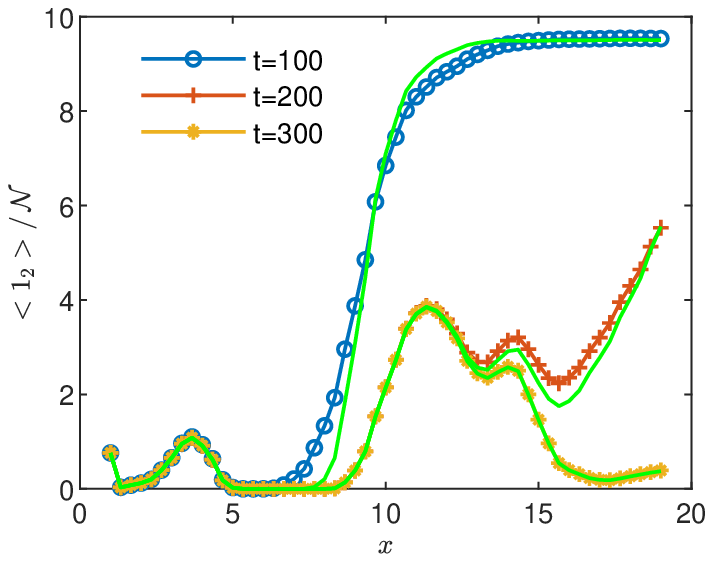}
\caption{\label{fig:aqv_cg2}Profiles of the CGST upscaled oxygen concentration at three sampling times compared to volume averages (full lines), sampled on the horizontal midline of the domain.}
\end{minipage}
\end{figure}

Figures~\ref{fig:aqv_c1}~and~\ref{fig:aqv_c2} show the meandering contours of the contaminant and oxygen concentrations transported in a sample of random velocity field. The CGST and the volume concentrations compared in Figs.~\ref{fig:aqv_cg1}~and~\ref{fig:aqv_cg2} indicate significant discrepancies negligible only for the oxygen concetrations. A more complete image of the considered simulation scenario is provided by the computed discrepancies presented in Table~\ref{table:Err_2Daqv}. Excepting the large maximum difference at $t=200$ (at very small concentration values), the discrepancy for the contaminant concentrations is at most a few percents. Significant discrepancies $e_{c_2}$ and $\varepsilon_{c_2}$, larger than 50\%, are instead observed for the averaged oxygen concentration.

\begin{table}
\caption{\label{table:Err_2Daqv}Relative differences between volume and CGST averages: two-dimensional saturated aquifer.}
\begin{center}
\begin{tabular}{ c | c c | c c }
   $ t $ & $e_{c_1}$ & $\varepsilon_{c_1}$ &  $e_{c_2}$ & $\varepsilon_{c_2}$ \\
\hline
     100 & 0.0423 & 0.1841 & 0.0544 & 0.6610 \\
     200 & 0.0010 & 0.9984 & 0.0802 & 0.2101 \\
     300 & 0.0035 & 0.0035 & 0.0040 & 0.0410 \\
  \hline
\end{tabular}
\end{center}
\end{table}

\section{Summary}
\label{sec:concl}

The CGST-GRW approach proposed in this article performs the space-time upscaling adapted to the spatio-temporal scale of a hypothetical measurement. The CGST averages over symmetrical volumes and time intervals of the microscopic description of a corpuscular system given by piecewise analytic time functions are a.e. continuous fields which verify identities similar to the local balance equations of the fluid dynamics. If the averages correspond to upscaled concentrations, these identities allow us to derive theoretical expressions for the flow velocity and the intrinsic diffusion coefficients of the system of particles. The microscopic description of reactive transport in porous media consists of GRW simulations of the advective-diffusive movement of the actual number of molecules involved in reactions. The chemical reactions are modeled deterministically by solving the reaction system with concentrations expressed as mole fractions of the number of reacting molecules. The main results of this study are summarized as follows.

\begin{itemize}
\item The codes for the CGST-GRW numerical modeling are verified by computing the theoretical expressions for flow velocity and intrinsic diffusion coefficients and one finds an excellent agreement with the nominal values of these input parameters.
\item If the reacting system consists of biomolecular reactions in a uniform flow, the initial distribution of the reacting species is uniform, and the boundary conditions force constant concentration solutions, the CGST average is practically indistinguishable from the classical volume average.
\item For biodegradation reactions in heterogeneous saturated aquifers, the CGST averages may significantly differ from the volume averages even if the flow is stationary, due to the time-dependence of the reaction process.
\item The discrepancy between volume and CGST averages may be dramatically increased if the biodegradation process takes place in partially saturated soils with time-dependent flow velocity.
\item According to Proposition~\ref{proposition:stoch_average} and Remark~\ref{rem:MCaverage}, the CGST average concentration is equivalent to a space-time average over the continuous field of the fine-grained concentration. This renders possible comparisons with time averages of classical volume sampling solutions obtained by integrating continuous fields or with time averages over solutions provided by homogenization methods. Kernel density estimates of concentrations in particle tracking simulations averaged over symmetrical time intervals can be compared as well with CGST averages computed with the same spatial and temporal scales.
\end{itemize}

\vspace{0.3cm}
\noindent{\bf\large Appendix A. CGST averaging on $d$-dimensional cubes}
\label{appA}
\vspace{0.3cm}

We consider the evolution of a system of $\mathcal{N}$ particles during the time interval $I=[0,T]$. Let $I_i=[t_{i}^{+},t_{i}^{-}]$, $0\leq t_{i}^{+}\leq t_{i}^{-}$, $i=1,\ldots,\mathcal{N}$, be the interval of existence of each particle, where $t_{i}^{+}$ and $t_{i}^{-}$ are the birth and dead times, respectively. These time points could be, for instance, those of generation and consumption of a molecule in chemical reactions. A {\it kinematic microscopic description} of a physical system can be achieved by arbitrary time functions $\varphi_i(t):I\mapsto\mathbb{R}$, with $\varphi_i(t)=0$ $\forall\; t\in I\setminus I_i$, $i=1,\ldots,\mathcal{N}$. Particular cases of functions $\varphi_i$ are the components of the particle's trajectory, $x_{\alpha i}:I\mapsto\mathbb{R}$, and the corresponding velocity components, $\xi_{\alpha i}:I\mapsto\mathbb{R}$, $\xi_{\alpha i}(t)=d x_{\alpha i}(t)/dt$. The existence intervals $I_i$ may contain finite number of points consisting of sets of zero Lebesgue measure $I_{i}^{*}=\{t_{i1},\ldots,t_{ik},t_{i(k+1)},\ldots t_{in}\}\subset I_i$ where the functions $\varphi_i$ undergo discontinuous variations. For instance, this is the case of velocity changes at successive time steps in random walk simulations. The restrictions of $\varphi_i$ to intervals $(t_{ik},t_{i(k+1)})$ are assumed to be analytic functions. So, we consider a kinematic description given by sets of piecewise analytic functions associated to the system of particles.

Let $C(\mathbf{x},a)=\{\mathbf{x},\mathbf{y}\in\mathbb{R}^d\; | \; |y_{\alpha}-x_{\alpha}|<a, \alpha=1,\ldots,d\}$ be an open $d$-dimensional cube with side $2a$, centered at $\mathbf{x}$. Similarly to the approach for coarse-grained averaging over phase points describing ensembles in statistical mechanics, for the case $I_i=I,\; \forall\; i=1,\ldots,\mathcal{N}$, presented in \cite[Chap. 5]{Vamos2007}, we define the CGST average over the cube $C(\mathbf{x},a)$ and the time interval $(t-\tau, t+\tau)$, $0<\tau<T/2$,
\begin{equation}\label{eq:A1}
\langle\varphi\rangle(\mathbf{x},t;a,\tau)=\frac{1}{2\tau(2a)^d}\sum_{i=1}^{\mathcal{N}}\int_{t-\tau}^{t+\tau}\varphi_i(t')\prod_{\alpha=1}^{d}
\chi_{\alpha}(x_{\alpha i}(t'))dt',
\end{equation}
where $\chi_{\alpha}$ is the characteristic function of the interval $(x_{\alpha}-a,x_{\alpha}+a)$.

A particle enters (leaves) the cube $C(\mathbf{x},a)$ at time instants $u_i$ which solve the system of equations
\begin{equation}\label{eq:A2}
x_{\alpha i}(u_i)-x_{\alpha} \pm a = 0,\; \alpha=1,2,\ldots,d.
\end{equation}
We denote by $U'_i=\{u'_{i1},u'_{i2},\ldots, u'_{in'}\}$ and $U''_i=\{u''_{i1},u''_{i2},\ldots, u''_{in''}\}$ the sets of time points corresponding to instances at which the $i$-th particles enters or leaves the cube $C(\mathbf{x},a)$. It follows that, for a fixed $\mathbf{x}$, the integrand in Eq.~(\ref{eq:A1}),
\begin{equation}\label{eq:A3}
G_i(\mathbf{x},t')=\varphi_i(t')\prod_{\alpha=1}^{d}
\chi_{\alpha}(x_{\alpha i}(t'))=\varphi_i(t')\chi_{C(\mathbf{x},a)}(\mathbf{x}_{i}(t')),
\end{equation}
is a continuous time function in $(t-\tau,t+\tau)$, except at a finite number of points $\{t_{i}^{+},t_{i}^{-}\}\bigcup U'_{i}\bigcup U''_{i}\bigcup I_{i}^{*}$ where it has finite jump discontinuities. Hence, $G(\mathbf{x},t')$ is a function with bounded variation and Riemann integrable.

Using Eq.~(\ref{eq:A1})~and~Eq.~(\ref{eq:A3}), the time derivative of the CGST average can be expressed as
\begin{equation}\label{eq:A4}
\partial_t \langle\varphi\rangle (\mathbf{x},t;a,\tau)=\frac{1}{2\tau(2a)^d}\sum_{i=1}^{\mathcal{N}}\left[G_i(\mathbf{x},t+\tau)-G_i(\mathbf{x},t-\tau)\right].
\end{equation}
The continuity of the time derivative (\ref{eq:A4}) is determined by that of the function $G_i$. According to (\ref{eq:A3}), $G_i$ is not continuous when $\varphi_i$ is discontinuous and $\chi_{C(\mathbf{x},a)}$ is nonvanishing or, conversely, $\chi_{C(\mathbf{x},a)}$ is discontinuous and $\varphi_i$ is nonvanishing. The first situation occurs if the particle is created or consumed at $t=t_{i}^{\pm}$ or undergoes discontinuous variations at times $t\in I^{*}_{i}$ inside the cube $C(\mathbf{x},a)$. The second situation occurs during the time intervals $\subset I_i\setminus I_{i}^{*}$ when the particle lies on the surface of the cube, $\partial C(\mathbf{x},a)$. Taking into account the time arguments of $G_i$ in (\ref{eq:A4}), $\partial_t \langle\varphi\rangle$ is discontinuous on the set $\Omega'=\cup_{i=1}^{\mathcal{N}}\Omega'_{i}$, where
\begin{align*}
\Omega'_{i}=
&\{\{t_{i}^{\pm}-\tau,t_{i}^{\pm}+\tau\}\cap(\tau,T-\tau)\} \times C(\mathbf{x}_{i}(t_{i}^{\pm}),a)\}\cup\\
&\{\{\{t_{ik}-\tau,t_{ik}+\tau\}\;|\;t_{ik}\in I_{i}^{*}\}\cap(\tau,T-\tau)\} \times C(\mathbf{x}_{i}(t_{i}^{*}),a)\}\cup\\
&\{\{[t_{i}^{+}\pm\tau,t_{i}^{-}\pm\tau]\cap(\tau,T-\tau)\} \times \partial C(\mathbf{x}_{i}(t\mp\tau),a)\}\cup\\
&\{\{\{[t_{i}^{+}\pm\tau,t_{i1}\pm\tau)\cup(t_{i1}\pm\tau,t_{i2}\pm\tau)\cup\ldots\\
&\cup(t_{in}\pm\tau,t_{i}^{-}\pm\tau]\}\cap(\tau,T-\tau)\} \times \partial C(\mathbf{x}_{i}(t\mp\tau),a)\}.
\end{align*}
Since $\Omega'$ has zero Lebesgue measure in $\mathbb{R}^d\times(\tau, T-\tau)$, the time derivative $\partial_t \langle\varphi\rangle$ is a.e. continuous.

As a function with bounded variation, $G_i$ may be decomposed as a sum of a continuous function and a jump function, $G_i=G'_i+G''_i$. Accordingly,
\begin{equation}\label{eq:A5}
\partial_t \langle\varphi\rangle=(\partial_t \langle\varphi\rangle)'+(\partial_t \langle\varphi\rangle)''.
\end{equation}
Since, except at discontinuity points, $G_i$ is an analytic function, according to the fundamental theorem of Lebesgue integral calculus, the continuous part $G'_i$ is absolutely continuous. Hence,
\begin{equation}\label{eq:A6}
(\partial_t \langle\varphi\rangle)'=\textstyle{\left\langle\frac{d\varphi}{dt}\right\rangle}.
\end{equation}
The second term of Eq.~(\ref{eq:A5}) consists of positive contributions when particles enter into the cube $C(\mathbf{x},a)$ and negative contributions when they leave the cube, as well as contributions of other three types of discontinuities, which will be described below. According to Eq.~(\ref{eq:A1}), $(\partial_t \langle\varphi\rangle)''$ can be written as
\begin{equation}\label{eq:A7}
(\partial_t \langle\varphi\rangle)''=\frac{1}{2\tau(2a)^d}\sum_{i=1}^{\mathcal{N}}\sum_{\alpha=1}^{d}\left[\sum_{u_{\alpha}\in W'_i}\varphi_i(u_{\alpha})-\sum_{u_{\alpha}\in W''_i}\varphi_i(u_{\alpha})\right]+(\partial_t \langle\varphi\rangle)_{gen}+(\partial_t \langle\varphi\rangle)_{disc}+\epsilon,
\end{equation}
where $W'_i=U'_i\bigcap(t-\tau,t+\tau)$ and $W''_i=U''_i\bigcap(t-\tau,t+\tau)$. The term $(\partial_t \langle\varphi\rangle)_{gen}$ accounts for generation/consumption of molecules in chemical reactions. It stems from variations of $G_i$ when the molecule is generated/consumed in the interior of $C(\mathbf{x},a)$ during the interval $(t-\tau,t+\tau)$:
\[
\Delta^{+}G_i=\varphi_i(t^+)\chi_{_{C(\mathbf{x},a)}}(\mathbf{x}_i(t^+))\chi_{_{(t-\tau,t+\tau)}}(t^+),\;\;
\Delta^{-}G_i=\varphi_i(t^-)\chi_{_{C(\mathbf{x},a)}}(\mathbf{x}_i(t^-))\chi_{_{(t-\tau,t+\tau)}}(t^-),
\]
where $\chi_{C(\mathbf{x},a)}$ and $\chi_{(t-\tau,t+\tau)}$ are the characteristic functions of the cube and of the time-averaging interval, respectively. Hence,
\begin{equation*}
(\partial_t \langle\varphi\rangle)_{gen}=\frac{1}{2\tau(2a)^d}\sum_{i=1}^{\mathcal{N}}\left[\Delta^{+}G_i+\Delta^{-}G_i\right].
\end{equation*}
The term $(\partial_t \langle\varphi\rangle)_{disc}$ describes variations due to the discontinuities of $\varphi_i$ at points from $I_{i}^*$ \cite{Vamosetal1997,Vamosetal2000}. This is, for instance, the case if $\varphi_i=\xi_{\alpha i}$ in random walk simulations when the velocity undergoes finite jumps at the end of each time step. In general,
\begin{equation*}
(\partial_t \langle\varphi\rangle)_{disc}=\frac{1}{2\tau(2a)^d}\sum_{i=1}^{\mathcal{N}}\sum_{s\in I_{i}^{*}}\left[\varphi_i(s-0)-\varphi_i(s+0)\right],
\end{equation*}
where $\varphi_i(s-0)$ and $\varphi_i(s+0)$ are the limits from the left and from the right respectively. The term $\epsilon$ cancels spurious contributions $\varphi_i(u_{\alpha})$ which occur when the trajectory intersects an vertex or an edge of the $d$-dimensional cube with $d\geq 2$. Since $\epsilon$ is defined on a set with zero Lebesgue measure, we assume that it is negligible. For instance in two-dimensional random walk simulations, unless in the worst case when most of the particles are located in the four corners of the square, the relative number of spurious contributions, of the order $1/\mathcal{L}$, where $\mathcal{L}$ is the number of lattice points per side, is  quite small and $\epsilon$  can be disregarded.

The CGST average defined by (\ref{eq:A1}) depends on $x_{\alpha}$ through the time points $u'_i$ and $u''_{i}$. The derivative of (\ref{eq:A2}) w.r.t. $x_{\alpha}$ leads to $\partial u_i / \partial x_{\alpha} = 1/(d x_{\alpha i}/du_i)=1/\xi_{\alpha i}$. If $u_i\in(t-\tau,t+\tau)$, then $u'_i$ is the lower integration limit of the time integral in (\ref{eq:A1}) and $u''_{i}$ is the upper limit. This implies the following expression of the partial space derivative,
\begin{equation}\label{eq:A10}
\partial_{x_{\alpha}} \langle\varphi\rangle=\frac{1}{2\tau(2a)^d}\sum_{i=1}^{\mathcal{N}}\left[-\sum_{u_{\alpha}\in W'_i}\frac{\varphi_i(u_{\alpha})}{\xi_{\alpha i}}+\sum_{u_{\alpha}\in W^{''}_i}\frac{\varphi_i(u_{\alpha})}{\xi_{\alpha i}}\right].
\end{equation}
As long as the particle's velocity $\xi_{\alpha i}\neq 0$ is defined and the implicit function theorem can be applied to define  $\partial u_i / \partial x_{\alpha} = 1/\xi_{\alpha i}$, the spatial derivative $\partial_{x_{\alpha}} \langle\varphi\rangle$ exists and is continuous. These conditions are not met on sets $\Omega''_{i}$ consisting of time intervals determined by $t\pm\tau$ with $t\in\{t_{i}^{+},t_{i}^{-}\}\cup I^{*}$ and sets $\{\mathbf{x}\in \partial C(\mathbf{x}_{i}(t),a)\}$ as well as in cases where particles move on the surface of the cube and Eq.~(\ref{eq:A2}) is identically satisfied. A detailed analysis and the construction of the sets $\Omega''_{i}$ is rather involved \cite[Chap. 4]{Vamos2007} but, in the present context, it is not necessary. Indeed, since the surface $\partial C(\mathbf{x}_{i}(t),a)$ has a null Lebesgue measure in $\mathbb{R}^d$, the union $\Omega''=\cup_{i=1}^{\mathcal{N}}\Omega''_{i}$ also has a null Lebesgue measure in $\mathbb{R}^d\times(\tau, T-\tau)$ and the spatial derivatives $\partial_{x_{\alpha}} \langle\varphi\rangle$ are a.e. continuous. The continuity of the partial derivatives $\partial_{x_{\alpha}} \langle\varphi\rangle$ and $\partial_{t} \langle\varphi\rangle$ ensures the a.e. continuity of the CGST averages $\langle\varphi\rangle$.

From (\ref{eq:A10}) and (\ref{eq:A7}) one obtains
\begin{equation}\label{eq:A11}
(\partial_t \langle\varphi\rangle)^{''}=-\sum_{\alpha=1}^{d}\partial_{x_{\alpha}} \langle\varphi\xi_{\alpha}\rangle+(\partial_t \langle\varphi\rangle)_{gen}+(\partial_t \langle\varphi\rangle)_{disc}.
\end{equation}
By inserting (\ref{eq:A6}) and (\ref{eq:A11}) into (\ref{eq:A5}), one obtains the following identity verified by the CGST averages,
\begin{equation}\label{eq:A12}
\partial_t \langle\varphi\rangle+\nabla\cdot\langle\varphi\boldsymbol{\xi}\rangle=\textstyle{\left\langle\frac{d\varphi}{dt}\right\rangle}+\delta\varphi,
\end{equation}
where $\delta\varphi=(\partial_t \langle\varphi\rangle)_{gen}+(\partial_t \langle\varphi\rangle)_{disc}$. Equation~(\ref{eq:A12}) is valid in $\mathbb{R}^d\times(\tau, T-\tau)\setminus(\Omega'\cup\Omega'')$, that is, a.e. in the space-time domain on which the CGST averages $\langle\varphi\rangle$ are defined.

\vspace{1cm}
\noindent{\bf\large Appendix B. Discrepancy between volume and CGST averages for different $a$ and $\tau$}
\label{appB}

\begin{figure}
\includegraphics[width=\linewidth]{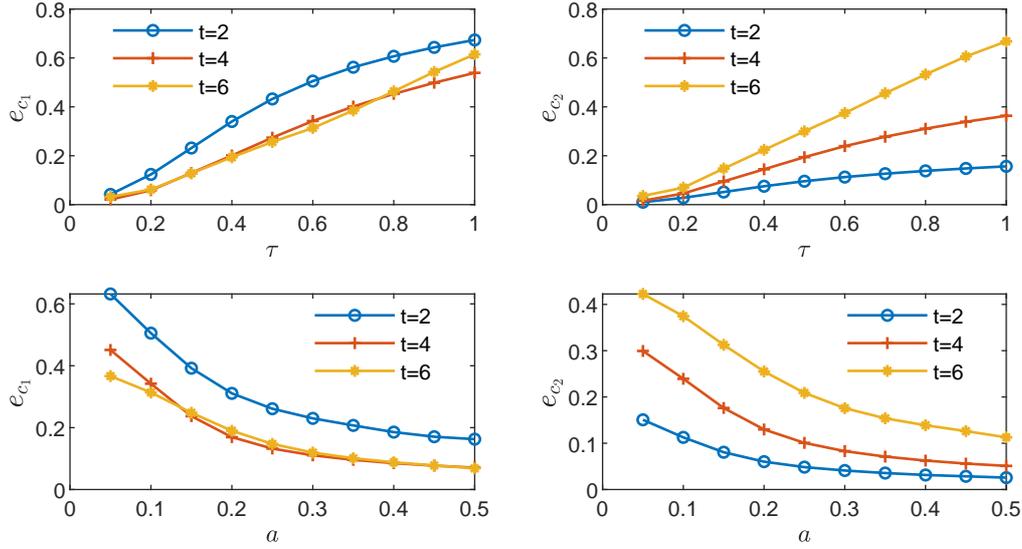}
\caption{\label{fig:errL2}Global relative differences $e_{c_{\nu}}$ between volume and CGST averages: increasing $\tau$, $a=0.1$ (top); increasing $a$, $\tau=0.6$ (bottom).}
\end{figure}

\begin{figure}
\includegraphics[width=\linewidth]{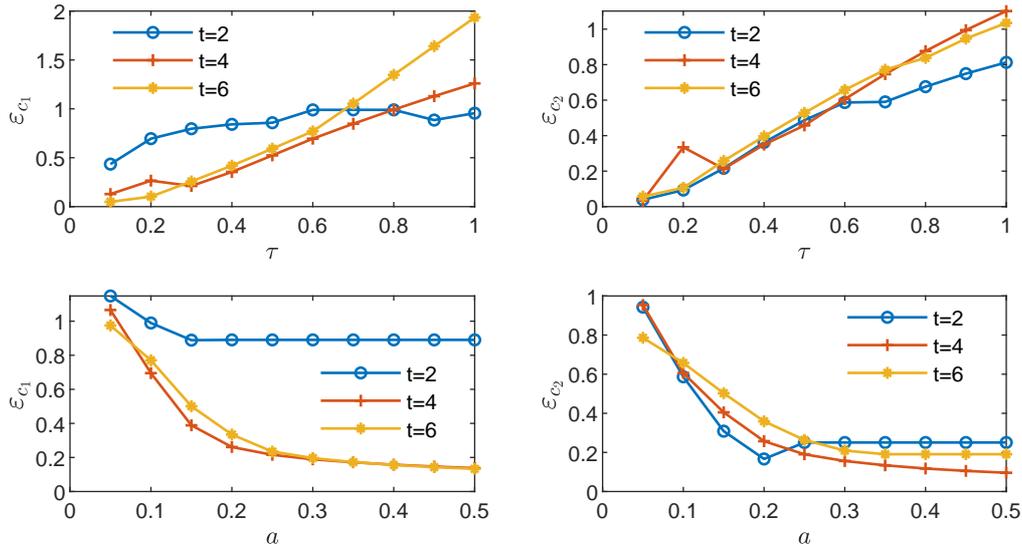}
\caption{\label{fig:err}Maximum discrepancy $\varepsilon_{c_{\nu}}$ between volume and CGST averages: increasing $\tau$, $a=0.1$ (top); increasing $a$, $\tau=0.6$ (bottom).}
\end{figure}

Discrepancies between the volume and CGST averages are computed for the one-dimensional example of reactive transport governed by Monod model presented in Section~\ref{sec:MonodReaction_soil}. The influence of the temporal scale $\tau$ and of the spatial scale $a$ on the global relative difference $e_{c_{\nu}}$ is illustrated in the sub-plots from Fig.~\ref{fig:errL2}. Similarly, the influence of $\tau$ and $a$ on $\varepsilon_{c_{\nu}}$ is illustrated in Fig.~\ref{fig:err}.

Figure Fig.~\ref{fig:errL2} shows a monotone increase of $e_{c_{\nu}}$ with $\tau$ for fixed $a$ and a monotone decrease with $a$ for fixed $\tau$. One remarks larger $e_{c_{\nu}}$ values for contaminant ($\nu=1$) and smaller values for oxygen ($\nu=2$) at the first sampling time (that is, close to the contaminant source)

Figure Fig.~\ref{fig:err} presents a rather irregular behavior of $\varepsilon_{c_{\nu}}$ at the first sampling time. The maximum difference is less influenced by the increase of $\tau$ and $a$. This is especially notable for the contaminant, with $\varepsilon_{c_{1}}$ values that remain almost constant with the decrease of $\tau$ and the increase of $a$.

\section*{Acknowledgements}
N. Suciu acknowledges the financial support of the Deutsche Forschungsgemeinschaft (DFG, German Research Foundation) under Grant SU 415/4-1 -- 405338726 ``Integrated global random walk model for reactive transport in groundwater adapted to measurement spatio-temporal scales''.



\begin{thebibliography}{00}
	
\bibitem{AkagiandOka2020}Akagi, G., Oka, T., 2020. Space-time homogenization for nonlinear diffusion. Preprint arXiv:2007.09977.
\href{https://doi.org/10.48550/arXiv.2007.09977}{https://doi.org/10.48550/arXiv.2007.09977}

\bibitem{Alecsaetal2020}Alecsa, C.D., Boros, I., Frank, F., Knabner, P., Nechita, M., Prechtel, A., Rupp, A., Suciu, N., 2019. Numerical benchmark study for fow in heterogeneous aquifers. Adv. Water Resour., 138, 103558.
\href{https://doi.org/10.1016/j.advwatres.2020.103558}{https://doi.org/10.1016/j.advwatres.2020.103558}

\bibitem{Andricevic1998}Andri\v{c}evi\'{c}, R., 1998. Effects of local dispersion and sampling volume on the evolution of concentration fluctuations in aquifers. Water Resour. Res., 34(5), 1115--29.
\href{http://dx.doi.org/10.1029/98WR00260}{http://dx.doi.org/10.1029/98WR00260}

\bibitem{AuriaultandAdler1995}Auriault, J.L., Adler, P.M., 1995. Taylor dispersion in porous media: analysis by multiple scale expansions. Advances in Water Resources, 18(4), 217-226.
\href{https://doi.org/10.1016/0309-1708(95)00011-7}{https://doi.org/10.1016/0309-1708(95)00011-7}

\bibitem{BarrosandFiori2014}de Barros, F.P.J., Fiori, A., 2014. First-order based cumulative distribution function for solute concentration in heterogeneous aquifers: Theoretical analysis and implications for human health risk assessment. Water Resour. Res., 50(5), 4018--4037.
\href{http://dx.doi.org/10.1002/2013WR015024}{http://dx.doi.org/10.1002/2013WR015024}

\bibitem{Bayer-Raichetal2004}Bayer-Raich, M., Jarsj\"{o}, J., Liedl, R., Ptak, T., Teutsch, G.2004. Average contaminant concentration and mass flow in aquifers from time-dependent pumping well data: Analytical framework. Water Resour. Res., 40(8), W08303.
\href{http://dx.doi.org/10.1029/2004WR003095}{http://dx.doi.org/10.1029/2004WR003095}

\bibitem{BauseandKnabner2004}Bause, M., Knabner, P., 2004. Numerical simulation of contaminant biodegradation by higher order methods and adaptive time stepping. Comput. Visual. Sci., 7(2), 61--78.
\href{https://doi.org/10.1007/s00791-004-0139-y}{https://doi.org/10.1007/s00791-004-0139-y}	

\bibitem{BensoussanLionsPapanicolaou2011}Bensoussan, A., Lions, J.L. and Papanicolaou, G., 2011. Asymptotic analysis for periodic structures (Vol. 374). American Mathematical Soc.

\bibitem{Berkowitz2021}Berkowitz, B., 2022. HESS Opinions: Chemical transport modeling in subsurface hydrological systems -- Space, time, and the ``holy grail'' of ``upscaling''. Hydrol. Earth. Syst. Sci., 26, 2161-2180.
\href{https://doi.org/10.5194/hess-26-2161-2022}{https://doi.org/10.5194/hess-26-2161-2022}

\bibitem{Brunneretal2012}Brunner, F., Radu, F.A., Bause, M., Knabner, P., 2012. Optimal order convergence
of a modified BDM1 mixed finite element scheme for reactive
transport in porous media. Adv. Water Resour., 35, 163--171.
\href{http://dx.doi.org/10.1016/j.advwatres.2011.10.001}{http://dx.doi.org/10.1016/j.advwatres.2011.10.001}

\bibitem{Bringedaletal2016}Bringedal, C., Berre, I., Pop, I. S., Radu, F. A., 2016. Upscaling of nonisothermal reactive porous media flow under dominant P\'{e}clet number: the effect of changing porosity. Multiscale Model. Simul. 14(1), 502--533.
\href{https://doi.org/10.1137/15M1022781}{https://doi.org/10.1137/15M1022781}

\bibitem{CaroniandFiorotto2005}Caroni, E., Fiorotto, V., 2005. Analysis of concentration as sampled in natural aquifers. Transport Porous Media, 59(1),19--45.
\href{http://dx.doi.org/10.1007/s11242-004-1119-x}{http://dx.doi.org/10.1007/s11242-004-1119-x}

\bibitem{Cirpkaetal1999}Cirpka, O.A., Frind, E.O.,  Helmig, R., 1999. Numerical simulation of biodegradation controlled by transverse mixing. J. Contam. Hydrol., 40(2), 159-182.
\href{https://doi.org/10.1016/S0169-7722(99)00044-3}{https://doi.org/10.1016/S0169-7722(99)00044-3}

\bibitem{Cushman1983}Cushman, J.H., 1983. Multiphase transport equations: I - general equation for macroscopic statistical, local space-time homogeneity. Transport Theor. Stat. Phys., 12(1), 35--71.
\href{http://dx.doi.org/10.1080/00411458308212731}{http://dx.doi.org/10.1080/00411458308212731}

\bibitem{DabrowskaRykala2021}Dabrowska, D., Rykala, W. (2021). A review of lysimeter experiments carried out on municipal landfill waste. Toxics, 9(2), 26.
\href{https://doi.org/10.3390/toxics9020026}{https://doi.org/10.3390/toxics9020026}

\bibitem{DestouniandGraham1997}Destouni, G., Graham, W., 1997. The influence of observation method on local concentration statistics in the subsurface. Water Resour .Res., 33(4),663--676.
\href{http://dx.doi.org/10.1029/96WR03955}{http://dx.doi.org/10.1029/96WR03955}

\bibitem{Eberhardetal2007}Eberhard, J.P., Suciu, N., Vamo\c{s}, C., 2007. On the self-averaging of dispersion for transport in quasi-periodic random media. J. Phys. Math. Theor., 40(4), 597--610.
\href{https://doi.org/10.1088/1751-8113/40/4/002}{https://doi.org/10.1088/1751-8113/40/4/002}

\bibitem{Fernandez-GarciaandSanchez-Vila2011}Fern\`{a}ndez-Garcia, D., Sanchez-Vila, X., 2011. Optimal reconstruction of concentrations, gradients and reaction rates from particle distributions. J. Contam. Hydrol., 120, 99--114.
\href{https://doi.org/10.1016/j.jconhyd.2010.05.001}{https://doi.org/10.1016/j.jconhyd.2010.05.001}

\bibitem{GrayandMiller2013}Gray, W.G., Miller, C.T., 2013. A generalization of averaging theorems for porous medium analysis. Adv. Water Resour. 62, 227--37.
\href{https://doi.org/10.1016/j.advwatres.2013.06.006}{https://doi.org/10.1016/j.advwatres.2013.06.006}

\bibitem{Grayetal2015}Gray, W.G., Dye, A.L., McClure, J.E., Pyrak-Nolte, L.J., Miller, C.T., 2015. On the dynamics and kinematics of two-fluid-phase flow in porous media. Water Resour. Res., 51(7), 5365-81.
\href{https://doi.org/10.1016/10.1002/2015WR016921}{https://doi.org/10.1016/10.1002/2015WR016921}

\bibitem{HeandSykes1996}He, Y., Sykes, J.F., 1996. On the spatial-temporal averaging method for modeling transport in porous media. Transport Porous Media, 22(1), 1--51.
\href{https://doi.org/10.1007/BF00974310}{https://doi.org/10.1007/BF00974310}

\bibitem{Hesseetal2009}He{\ss}e, F., Radu, F.A., Thullner, M. and Attinger, S., 2009. Upscaling of the advection-diffusion-reaction equation with Monod reaction. Adv. Water Resour., 32(8), 1336--1351.
\href{https://doi.org/10.1016/j.advwatres.2009.05.009}{https://doi.org/10.1016/j.advwatres.2009.05.009}

\bibitem{Hornung1996}Hornung, U., 1996. Homogenization and porous media (Vol. 6). Springer Science \& Business Media.
\href{https://link.springer.com/book/10.1007/978-1-4612-1920-0}{https://link.springer.com/book/10.1007/978-1-4612-1920-0}

\bibitem{Illiano2020}Illiano, D., Pop, I.S., Radu, F.A., 2020. Iterative schemes for surfactant transport in porous media. Comput. Geosci.
\href{https://doi.org/10.1007/s10596-020-09949-2}{https://doi.org/10.1007/s10596-020-09949-2}

\bibitem{IrvingandKirkwood1950}Irving, J.H., Kirkwood, J.G., 1950. The statistical mechanical theory of transport processes. IV. The equations of hydrodynamics. J. Chem. Phys., 18(6), 817--829.
\href{https://doi.org/10.1063/1.1747782}{https://doi.org/10.1063/1.1747782}

\bibitem{Klofkornetal2002}Kl\"{o}fkorn, R., Kr\"{o}ner, D., Ohlberger, M., 2002. Local adaptive methods for convection dominated problems. Int. J. Numer. Meth. Fluid., 40(1-2), 79--91.
\href{https://doi.org/10.1002/fld.268}{https://doi.org/10.1002/fld.268}

\bibitem{Kumaretal2016}Kumar, K., Neuss-Radu, M., Pop, I.S., 2016. Homogenization of a pore scale model for precipitation and dissolution in porous media. IMA J. Appl. Math., 81(5), 877--897.
\href{https://doi.org/10.1093/imamat/hxw039}{https://doi.org/10.1093/imamat/hxw039}

\bibitem{ListandRadu2016} List, F. and Radu, F.A., 2016. A study on iterative methods for solving Richards' equation. Comput. Geosci.  20 (2), 341--353.
\href{https://doi.org/10.1007/s10596-016-9566-3}{https://doi.org/10.1007/s10596-016-9566-3}

\bibitem{Maillouxetal2003}Mailloux, B.J., Fuller, M.E., Rose, G.F., Onstott, T.C., DeFlaun, M.F., Alvarez, E.,
Hemingway, C., Hallet, R.B., Phelps, T.J., Griffin, T., 2003. A modular injection system, multilevel sampler, and manifold for tracer tests. Ground
Water, 41(6), 816--27.
\href{http://dx.doi.org/10.1111/j.1745-6584.2003.tb0242}{http://dx.doi.org/10.1111/j.1745-6584.2003.tb0242}

\bibitem{McClure2017}McClure, J.E., Dye, A.L., Miller, C.T., Gray, W.G., 2017. On the consistency of scale among experiments, theory, and simulation. Hydrol. Earth Syst. Sci., 21(2), 1063.
\href{http://dx.doi.org/10.5194/hess-21-1063-2017}{10.5194/hess-21-1063-2017}

\bibitem{Nolenetal2008}Nolen, J., Papanicolaou, G., Pironneau, O., 2008. A framework for adaptive multiscale methods for elliptic problems. Multiscale Model. Simul., 7, 171--196.
\href{http://dx.doi.org/10.1137/070693230}{http://dx.doi.org/10.1137/070693230}

\bibitem{Nordbottenetal2007}Nordbotten, J.M., Celia, M.A., Dahle, H.K., Hassanizadeh, S.M., 2007. Interpretation of macroscale variables in Darcy's law. Water Resour. Res., 43(8), W08430.
\href{http://dx.doi.org/10.1029/2006WR005018}{http://dx.doi.org/10.1029/2006WR005018}

\bibitem{Popetal2004} Pop, I.S., Radu, F.A., Knabner, P., 2004. Mixed finite elements for the Richards' equation: linearization procedure. J. Comput. Appl. Math., 168 (1-2), 365--373.
\href{https://doi.org/10.1016/j.cam.2003.04.008}{https://doi.org/10.1016/j.cam.2003.04.008}

\bibitem{Portaetal2012}Porta, G.M., Riva, M., Guadagnini, A., 2012. Upscaling solute transport in porous media in the presence of an irreversible bimolecular reaction. Adv. Water Resour., 35, 151--162.
\href{https://doi.org/10.1016/j.advwatres.2011.09.004}{https://doi.org/10.1016/j.advwatres.2011.09.004}

\bibitem{Schueleretal2016}Sch\"{u}ler, L., Suciu, N., Knabner, P., Attinger, S. 2016. A time dependent mixing model to close PDF equations for transport in heterogeneous aquifers. Adv. Water Resour., 96, 55--67.
\href{http://dx.doi.org/10.1016/j.advwatres.2016.06.012}{http://dx.doi.org/10.1016/j.advwatres.2016.06.012}

\bibitem{Schwarzeetal2001}Schwarze, H., Jaekel, U., Vereecken, H., 2001. Estimation of macrodispersion by different approximation methods for flow and transport in randomly heterogeneous media. Transp. Porous Media, 43(2), 265-287.
\href{https://doi.org/10.1023/A:1010771123844}{https://doi.org/10.1023/A:1010771123844}

\bibitem{Schwedeetal2008}Schwede, R.L., Cirpka, O.A., Nowak, W., Neuweiler, I., 2008. Impact of sampling volume on the probability density function of steady state concentration. Water Resour. Res., 44(12), W12433.
\href{http://dx.doi.org/10.1029/2007WR006668}{http://dx.doi.org/10.1029/2007WR006668}

\bibitem{Skoienetal2003}Sk{\o}ien, J.O., Bl\"{o}schl, G., Western, A.W., 2003. Characteristic space scales and timescales in hydrology. Water Resour. Res., 39(10).
\href{http://dx.doi.org/10.1029/2002WR001736}{http://dx.doi.org/10.1029/2002WR001736}

\bibitem{Srzicetal2013}Srzic, V., Cvetkovic, V., Andricevic, R., Gotovac, H., 2013. Impact of aquifer heterogeneity structure and local-scale dispersion on solute concentration uncertainty. Water Resour. Res., 49(6), 3712--3728.
\href{http://dx.doi.org/10.1002/wrcr.20314}{http://dx.doi.org/10.1002/wrcr.20314}

\bibitem{Sole-MariandFernandez-Garcia2018}Sole-Mari, G., Fern\`{a}ndez-Garcia, D., 2018, Lagrangian modeling of reactive transport in heterogeneous porous media with an automatic locally adaptive particle support volume, Water Resour. Res. 54(10), 8309--8331,
\href{http://dx.doi.org/10.1029/2018WR023033}{http://dx.doi.org/10.1029/2018WR023033}

\bibitem{Sole-Marietal2019}Sole-Mari, G., Bolster, D., Fern\`{a}ndez-Garcia, D., Sanchez-Vila, X., 2019. Particle density estimation with grid-projected and boundary-corrected adaptive kernels. Adv. Water Resour. 131, 103382.
\href{https://doi.org/10.1016/j.advwatres.2019.103382}{https://doi.org/10.1016/j.advwatres.2019.103382}

\bibitem{Suciu2001}Suciu, N.,  2001. On the connection between the microscopic and macoscopic modeling of the thermodynamic processes. Ed. Univ.
Pite\c{s}ti (in Romanian). English abstract: \href{http://ictp.acad.ro/suciu/abstract-PhD-Thesis-Suciu.pdf}{http://ictp.acad.ro/suciu/abstract-PhD-Thesis-Suciu.pdf}.

\bibitem{Suciuetal2015}Suciu, N., Radu, F.A., Attinger, S., Sch\"{u}ler, L., Knabner, P., 2015. A Fokker-Planck approach for probability distributions of species concentrations transported in heterogeneous media. J. Comput. Appl. Math., 289, 241--252.
\href{http://dx.doi.org/10.1016/j.cam.2015.01.030}{http://dx.doi.org/10.1016/j.cam.2015.01.030}

\bibitem{Suciuetal2016}Suciu, N., Sch\"{u}ler, L., Attinger, S., Knabner, P., 2016. Towards a filtered density function approach for reactive transport in groundwater. Adv. Water Resour., 90, 83--98.
\href{http://dx.doi.org/10.1016/j.advwatres.2016.02.016}{http://dx.doi.org/10.1016/j.advwatres.2016.02.016}

\bibitem{Suciu2019}Suciu, N., 2019. Diffusion in Random Fields. Applications to Transport in Groundwater. Birkh\"{a}user, Cham.
\href{https://doi.org/10.1007/978-3-030-15081-5}{https://doi.org/10.1007/978-3-030-15081-5}

\bibitem{Suciuetal2021}Suciu, N., Illiano, D., Prechtel, A., Radu, F.A., 2021. Global random walk solvers for fully coupled flow and transport in saturated/unsaturated porous media.  Adv. Water Resour., 152, 103935.
\href{https://doi.org/10.1016/j.advwatres.2021.103935}{https://doi.org/10.1016/j.advwatres.2021.103935}

\bibitem{SuciuandRadu2021}Suciu, N., Radu, F.A., 2022. Global random walk solvers for reactive transport and biodegradation processes in heterogeneous porous media. Adv. Water Resour. 166, 104268.
    \href{http://dx.doi.org/10.1016/j.advwatres.2022.104268}{http://dx.doi.org/10.1016/j.advwatres.2022.104268}

\bibitem{Suciuetal2023}Suciu, N., Radu, F.A., Pop, I.S., 2023. \href{https://github.com/PMFlow/SpaceTimeUpscaling}{https://github.com/PMFlow/SpaceTimeUpscaling}. Git repository. 

\bibitem{Teutschetal2000}Teutsch, G., Ptak, T., Schwarz, R., Holder, T., 2000. Ein neues Verfahren zur Quantifizierung der Grundwasserimmission: I. Theoretische Grundlagen. Grundwasser, 4(5), 170--175.
\href{http://dx.doi.org/10.1007/s767-000-8368-7}{http://dx.doi.org/10.1007/s767-000-8368-7}

\bibitem{ToninaandBelin2008}Tonina, D., Bellin, A., 2008. Effects of pore-scale dispersion, degree of heterogeneity, sampling size, and source volume on the concentration moments of conservative solutes in heterogeneous formations. Adv. Water Resour., 31(2), 339--354.
\href{http://dx.doi.org/doi:10.1016/j.advwatres.2007.08.009}{http://dx.doi.org/doi:10.1016/j.advwatres.2007.08.009}

\bibitem{Vamosetal1996}Vamo\c{s} C., Georgescu A., Suciu N., Turcu I., 1996. Balance equations for physical systems with corpuscular structure.
Phys. Stat. Mech. Appl., 227(1), 81--92.
\href{http://dx.doi.org/10.1016/0378-4371(95)00373-8}{http://dx.doi.org/10.1016/0378-4371(95)00373-8}.

\bibitem{Vamosetal1997}Vamo\c{s} C., Suciu N,, Georgescu, A., 1997. Hydrodynamic equations for
one-dimensional systems of inelastic particles. Phys. Rev. E, 55(5), 6277--6280.
\href{http://doi.org/10.1103/PhysRevE.55.6277}{http://doi.org/10.1103/PhysRevE.55.6277}


\bibitem{Vamosetal2000}Vamo\c{s}, C., Suciu, N., Blaj, W., 2000. Derivation of one-dimensional hydrodynamic model for stock price evolution.
Phys. Stat. Mech. Appl., 287(3), 461--467.
\href{http://dx.doi.org/10.1016/S0378-4371(00)00385-X}{http://dx.doi.org/10.1016/S0378-4371(00)00385-X}

\bibitem{Vamosetal2003}Vamo\c{s}, C., Suciu, N., Vereecken, H., 2003. Generalized random walk algorithm for the numerical modeling of complex diffusion processes. J. Comput. Phys.,  186, 527--544.
\href{https://doi.org/10.1016/S0021-9991(03)00073-1}{https://doi.org/10.1016/S0021-9991(03)00073-1}

\bibitem{Vamos2007}Vamo\c{s} C., 2007. Continuous fields associated to corpuscular systems. Risoprint, Cluj-Napoca (in Romanian).
English version: \href{http://ictp.acad.ro/vamos/cv-thesis/}{http://ictp.acad.ro/vamos/cv-thesis/}

\bibitem{VanDuijnandvanderZee2018}Van Duijn, C.J., van der Zee, S.E.A.T.M., 2018. Large time behaviour of oscillatory nonlinear solute transport in porous media. Chemical Engineering Science, 183, 86-94.
\href{https://doi.org/10.1016/j.ces.2018.02.045}{https://doi.org/10.1016/j.ces.2018.02.045}

\bibitem{Wiedemeieretal1999}Wiedemeier, T.H., Rifai, H.S., Newell, C.J., Wilson, J.T., 1999. Natural attenuation of fuels and chlorinated solvents in the subsurface. John Wiley \& Sons.

\bibitem{Wrightetal2021}Wright, E. E., Sund, N. L., Richter, D. H., Porta, G. M., Bolster, D., 2021. Upscaling bimolecular reactive transport in highly heterogeneous porous media with the LAgrangian Transport Eulerian Reaction Spatial (LATERS) Markov model. Stoch. Environ. Res. Risk. Assess., 35, 1529--1547.
\href{https://doi.org/10.1007/s00477-021-02006-z}{https://doi.org/10.1007/s00477-021-02006-z}

\end{thebibliography}
\end{document}